\documentclass{article}

\usepackage{microtype}
\usepackage{graphicx}
\usepackage{subfigure}
\usepackage{booktabs} 

\usepackage{hyperref}
\usepackage{url}
\usepackage{graphicx}
\usepackage{caption}
\usepackage{subcaption}
\usepackage{tikz}
\usetikzlibrary{positioning}
\usetikzlibrary{arrows.meta}

\newcommand{\vertiii}[1]{{\left\vert\kern-0.25ex\left\vert\kern-0.25ex\left\vert #1
\right\vert\kern-0.25ex\right\vert\kern-0.25ex\right\vert}}

\usepackage[accepted]{icml2024}

\usepackage{amsmath}
\usepackage{amssymb}
\usepackage{mathtools}
\usepackage{amsthm}

\usepackage[capitalize,noabbrev]{cleveref}

\theoremstyle{plain}

\theoremstyle{definition}

\theoremstyle{remark}

\newtheorem*{astral*}{Astral loss}

\usepackage[textsize=tiny]{todonotes}

\icmltitlerunning{Neural functional a posteriori error estimates}

\begin{document}

\twocolumn[
\icmltitle{Neural functional a posteriori error estimates}



\icmlsetsymbol{equal}{*}
\begin{icmlauthorlist}
\icmlauthor{Vladimir Fanaskov}{yyy}
\icmlauthor{Alexander Rudikov}{yyy}
\icmlauthor{Ivan Oseledets}{yyy,xxx}

\end{icmlauthorlist}

\icmlaffiliation{yyy}{Computational Intelligence Research Group, Center for Artificial Intelligence Technology, Skoltech, Moscow, Russia 121205}
\icmlaffiliation{xxx}{Artificial Intelligence
Research Institute}

\icmlcorrespondingauthor{Vladimir Fanaskov}{v.fanaskov@skoltech.ru}

\icmlkeywords{Neural Operators, PDEs, scientific machine learning, a posteriori error analysis, physics-informed machine learning}

\vskip 0.3in
]



\printAffiliationsAndNotice{\icmlEqualContribution} 

\begin{abstract}
We propose a new loss function for supervised and physics-informed training of neural networks and operators that incorporates a posteriori error estimate. More specifically, during the training stage, the neural network learns additional physical fields that lead to rigorous error majorants after a computationally cheap postprocessing stage. Theoretical results are based upon the theory of functional a posteriori error estimates, which allows for the systematic construction of such loss functions for a diverse class of practically relevant partial differential equations. From the numerical side, we demonstrate on a series of elliptic problems that for a variety of architectures and approaches (physics-informed neural networks, physics-informed neural operators, neural operators, and classical architectures in the regression and physics-informed settings), we can reach better or comparable accuracy and in addition to that cheaply recover high-quality upper bounds on the error after training.
\end{abstract}

\section{Introduction}
Recently it has been a surge of interest in cheap surrogates for partial-differential equation (PDEs) solvers in the machine learning (ML) community. This interest led to the construction of a plethora of novel efficient architectures including Fourier Neural Operator (FNO) \cite{li2020fourier}, Deep Operator Network (DeelONet) \cite{lu2019deeponet}, wavelet- and transformer-based approaches \cite{gupta2021multiwavelet}, \cite{hao2023gnot}.

Arguably, for low-dimensional case the most promising results appeared in the regression setting, when the neural network learns a mapping from the input PDEs data to the quantity of interest (e.g., the solution at a particular time or on the whole spacetime interval of interest), given observed data or results of the classical simulations. Often, it is possible to retain enough accuracy (the typical error is a few percent of the relative $L_2$ norm) and simultaneously reduce solution time by a few orders of magnitude \cite{li2022fourier}, \cite{brandstetter2022message}, \cite{lam2022graphcast}. In this scenario, one can hope to use neural network-based PDE solvers in compute-intensive applications such as weather forecasts, PDE-constraint optimization \cite{biegler2003large}, or inverse problems \cite{tarantola2005inverse}.

Despite being promising, neural PDE solvers are unreliable. Even when the results on approximation capabilities are available \cite{kovachki2021universal}, \cite{lanthaler2022error}, they do not guarantee that it is possible to reach good accuracy with standard training approaches (see \cite{fokina2023growing}, \cite{colbrook2022difficulty}). Since a priori error analysis seems futile, the natural alternative is a posteriori error analysis.\footnote{We briefly recall the difference between a priori and a posteriori error analysis. The former one is an estimation that does not depend on the obtained solution, e.g., it can be a statement on the order of convergence: $\text{error} \leq C h^p$, where $h$ is a grid spacing, $p\in\mathbb{N}$, and $C$ is an unspecified constant that depends on the problem data. The latter one is a concrete estimation that explicitly depends on the obtained approximate solution but not on the exact solution, i.e., $\text{error} \simeq f(\text{approximate solution}, \text{problem data})$.} If reliable and computationally cheap a posteriori error analysis were possible, it would allow the safe application of neural PDE solvers.

With this note, we hope to draw the attention of the ML community to a well-developed powerful approach to a posteriori error analysis known as \textit{functional a posteriori error analysis} \cite{neittaanmaki2004reliable}, \cite{mali2013accuracy}, \cite{repin2008posteriori}, \cite{Muzalevsky2021aposteriori}. The main strength of the technique is that it allows for error estimation regardless of the method used to construct approximation (contrast this with highly specialized approaches, e.g., FEM a posteriori analysis \cite{ainsworth1997posteriori}). This feature makes functional a posteriori analysis especially appealing for solvers based on neural networks that often act as black boxes.

More specifically, functional a posteriori error analysis supplies a technique to construct error majorant (upper bound on the error) that depends on approximate solution, PDE data, and additional fields that can be used to tighten the upper bound. The main properties of the functional approach are that: (i) the majorant remains upper bound for arbitrary approximate solution\footnote{The only requirement is that approximate solution should lay in the appropriate PDE-dependent functional space.} regardless of the quality and the nature of the approximation, (ii) it provides a tight upper bound which is saturated for exact solution only. Using this powerful approach, we contribute the following:
\begin{enumerate}
	\item In \cref{Section: Functional a posteriori error estimate as a novel loss function}, following pioneering contribution \cite{Muzalevsky2021aposteriori}, we define novel loss function for physics-informed training called Astral (neur\textbf{A}l a po\textbf{ST}erio\textbf{R}i function\textbf{A}l \textbf{L}oss).
	\item In \cref{Section: Comparison of losses for physics-informed neural networks}, we demonstrate that for elliptic equations training with Astral loss is much more robust than training with residual and variational losses (see \cref{fig:PiNN losses}). Besides that, Astral is equivariant to rescaling and allows for direct error control, which is impossible with other approaches.
	\item As we explain in \cref{Section: Functional a posteriori error estimate as a novel loss function} Astral can also be used for parametric PDEs. We present two schemes to achieve that: the physics-informed (unsupervised) \cref{fig:unsupervised UQ training} and supervised \cref{fig:supervised UQ training}. We test these schemes in \cref{Section: Application to parametric PDEs} and find that in unsupervised setting Astral outperforms residual training used in PINO \cite{li2021physics} by a significant margin (see \cref{fig:PINO vs unsupervised upper bound}).
\end{enumerate}

\section{Functional a posteriori error estimate as a novel loss function}
\label{Section: Functional a posteriori error estimate as a novel loss function}

\subsection{High-level description of the error estimate}
\label{Section: Functional a posteriori error estimate as a novel loss function. SUbsection: High-level description of the error estimate.}
Consider PDE in the abstract form
\begin{equation}
    \label{eq:abstract_PDE}
    \mathcal{A}\left[u, \mathcal{D}\right] = 0,
\end{equation}
where $\mathcal{A}$ is a nonlinear operator containing partial derivatives of the solution $u$, $\mathcal{D}$ stands for supplementary data such as initial conditions, boundary conditions, and PDE parameters.

Typically, the exact solution to \cref{eq:abstract_PDE} is not available, so one obtains only approximate solution $\widetilde{u}$. Now, it is desirable to estimate the quality of $\widetilde{u}$, i.e., the distance between approximate and exact solution in some norm $\left\| \widetilde{u}-u_{\text{exact}}\right\|$. The estimate (i) should have definitive relation to the error norm (e.g., upper or lower bound), (ii) is computable only from PDE data and approximate solution $\widetilde{u}$, (iii) is cheaper to compute than the approximate solution $\widetilde{u}$ itself. Estimates following desiderata are known for specific discretizations, most notably for finite element methods  \cite{babuvska1978posteriori}, but they are not applicable when the solution comes from neural networks. The only approach suitable for deep learning is discretization-agnostic functional a posteriori error analysis, which we now outline.

In general, functional a posteriori error bounds for a given PDE have a form
\begin{equation}
    \label{eq:abstract_functional_error_estimate}
    L[\widetilde{u}, \mathcal{D}, w_{L}]\leq \left\|\widetilde{u} -u_{\text{exact}}\right\|\leq U[\widetilde{u}, \mathcal{D}, w_{U}],
\end{equation}
where $u$ is an approximate solution, $w_{L}$ and $w_{U}$ are arbitrary free functions (we call them \textit{certificates} in the text below) from certain problem-dependent functional space, $\mathcal{D}$ is problem data (e.g., diffusion coefficient, viscosity, information on the geometry of the domain e.t.c.), $U$ and $L$, error majorant and minorant, are problem-dependent nonlinear functionals (see \cref{Section: Functional a posteriori error estimate as a novel loss function. Subsection: Concrete error bounds for elliptic PDEs.} for concrete example) \cite{mali2013accuracy}. Majorant and minorant has special properties:
\begin{enumerate}
    \item They are defined in a continuous sense for arbitrary $\widetilde{u}$, $w_{U}$ from certain functional spaces, that is, they do not contain information on a particular solution method, grid quantities, convergence, or smoothness properties.
    \item The bounds are tight. That is, solving $\sup_{w_{L}, u} L[u, \mathcal{D}, w_{L}] \text{ or } \inf_{w_{U}, u} U[u, \mathcal{D}, w_{U}]$ one recovers exact solution $u_{\text{exact}}$.
    \item They are explicitly computable based on $\widetilde{u}$ and the problem data $D$.
\end{enumerate}

Typical application of \cref{eq:abstract_functional_error_estimate} is to plug given approximate solution $\widetilde{u}$ and optimize for additional field $
    \sup_{w_{L}}L[\widetilde{u}, \mathcal{D}, w_{L}]\leq \left\|\widetilde{u}-u_{\text{exact}}\right\|\leq \inf_{w_{U}}U[\widetilde{u}, \mathcal{D}, w_{U}],$
to find error bounds as tight as possible. The optimization problems are potentially computationally demanding, so one can use deep learning to mitigate numerical costs.

The other possibility is to consider the upper bound in \cref{eq:abstract_functional_error_estimate} as a novel loss function and optimize it for $u$ and $w_{U}$. Since the bounds are tight, $u$ converges to the exact solution, and since the majorant provides an upper bound, we have direct access to the upper bound on error during training.

We provide concrete examples of error bounds for elliptic equations and discuss how to construct them for other equations in \cref{Section: Functional a posteriori error estimate as a novel loss function. Subsection: Concrete error bounds for elliptic PDEs.}, but before that we digress to explain in more details two approaches to combine deep learning with a posteriori functional error estimate outlined above.

\subsection{From error majorant to deep learning}
\label{Section: Functional a posteriori error estimate as a novel loss function. Subsection: How to incorporate error estimates into deep-learning pipelines.}

\begin{figure*}
  \centering
  \input{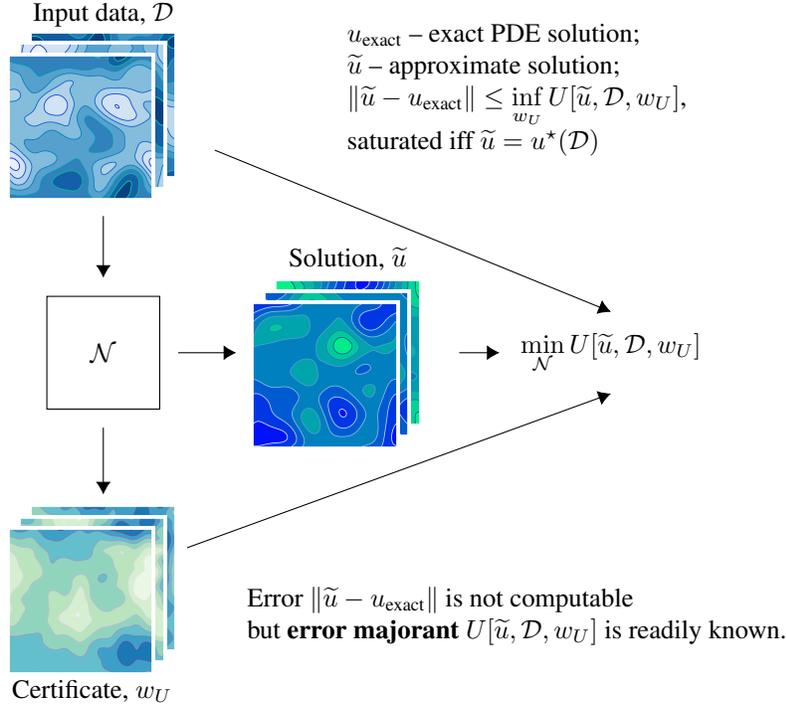}
  \caption{Unsupervised learning scheme with a posteriori functional error estimate and the explanation of basic properties of the upper bound. Neural network $\mathcal{N}$ takes input PDE data $\mathcal{D}$ and outputs approximate solution $\widetilde{u}$ and certificate $w_{U}$. Input to $\mathcal{N}$ and both outputs are plugged into \textbf{error majorant} $U[u, \mathcal{D}, w_{U}]$ which is used as an unsupervised loss.}
  \label{fig:unsupervised UQ training}
\end{figure*}

As we mentioned in the introduction, one of the most successful applications of deep learning to PDEs is the construction of cheap surrogates for classical solvers, i.e., the neural network produces an approximate solution given PDE data as input. To train such a network, one collects a dataset with pairs $(\text{PDE data}, \text{solution})$ using a classical solver and trains a neural network in the regression setting. One downside of this approach is the absence of guarantees. When neural networks provide solutions for unseen data, how can one measure the approximation error?

One possibility is to plug the solution $\widetilde{u}$ in discretization-agnostic upper bound $U[\widetilde{u}, \mathcal{D}, w_{U}]$ and optimize for $w_{U}$. This was done in \cite{Muzalevsky2021aposteriori} for solution obtained with physics-informed neural networks. However, the complexity of the optimization problem can be comparable with solving the original PDE, so this method does not look computationally appealing.

The natural alternative is find the solution $\widetilde{u}$ and the certificated $w_{U}$ in a single optimization run. For that we define novel loss function and the training strategy.
\begin{astral*}
    For PDE with functional error majorant $U[\widetilde{u}, \mathcal{D}, w_{U}]$, and neural network $\mathcal{N}(\mathcal{D}, \theta) = \left(\widetilde{u}, w_{U}\right)$ with weights $\theta$ that predicts solution $\widetilde{u}$ and certificate $w_{U}$ optimize
    \begin{equation}
        \label{eq:ASTRAL}
        \min_{\theta} U[\widetilde{u}, \mathcal{D}, w_{U}] \text{ s.t. } \left(\widetilde{u}, w_{U}\right) = \mathcal{N}(\mathcal{D}, \theta).
    \end{equation}
\end{astral*}
That is, training is done the same way as for classical physics-informed neural networks and neural operators, but with two differences. First, the loss function is the error majorant, not the residual. Second, the neural network has additional outputs since it needs to predict certificates $w_{U}$ along with the approximate solution $\widetilde{u}$. The concrete example of Astral loss is given in \cref{Section: Functional a posteriori error estimate as a novel loss function. Subsection: Concrete error bounds for elliptic PDEs.} for elliptic problem (see \cref{eq:upper_bound_nD}), it is also explained there how to construct it for other PDEs.

\begin{figure*}
  \centering
  \input{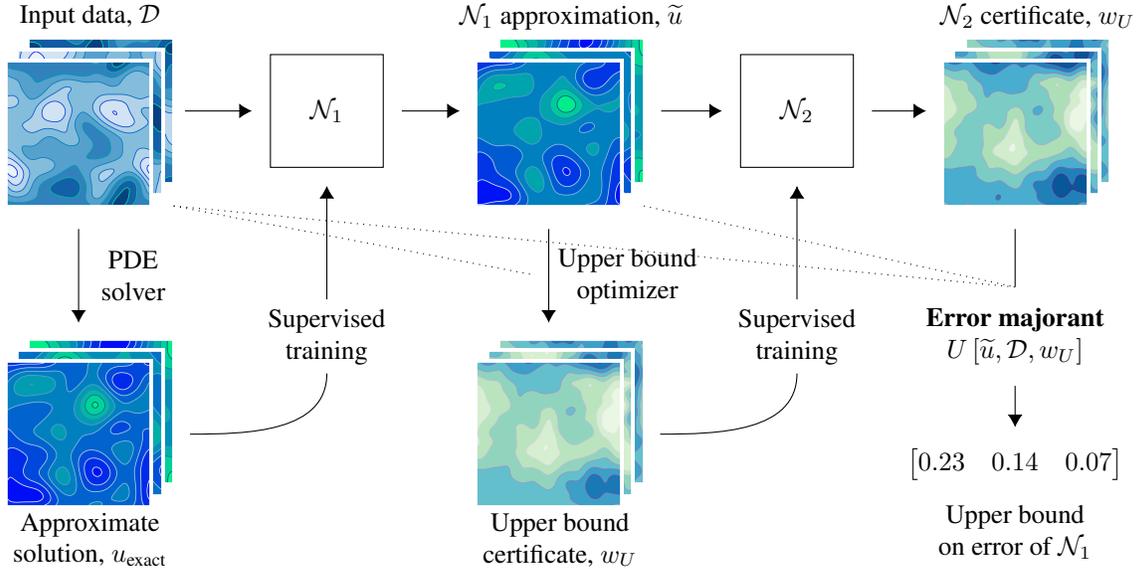}
  \caption{Supervised learning scheme with a posteriori functional error estimate. From left to right: (1) PDE is solved for selected input data; (2) pairs (input data, approximate solution) are used to train $\mathcal{N}_1$ with ordinary $L_2$ loss; (3) approximate solutions produced by $\mathcal{N}_1$ are plugged into upper bound and optimized certificates are found with direct optimization; (4) triples (input data, $\mathcal{N}_1$ approximation, certificate) are used to train $\mathcal{N}_2$ with ordinary $L_2$ loss; (5) certificates produced by $\mathcal{N}_2$, output of $\mathcal{N}_1$ and input data are plugged into upper bound to obtain bounds on error.}
  \label{fig:supervised UQ training}
\end{figure*}

Unfortunately, physics-informed approach are often computationally less efficient for low dimensional problems \cite{grossmann2023can}, \cite{karnakov2022optimizing}. So one can resort to standard solvers and optimizers and simply replace the computationally demanding search for optimal $w_{U}$ with a cheap neural network surrogate. This may sound suspicious. One can argue that we will lose guarantees on the upper bound if we outsource optimization to possibly inaccurate neural networks. Fortunately, this is not the case since the upper bound holds for arbitrary $w_{U}$. When neural networks fail to predict sufficiently accurate $w_{U}$, the upper bound merely ceases to be tight but remains a genuine upper bound. These considerations lead to the approach summarized in \cref{fig:supervised UQ training}. It consists of four stages:
\begin{enumerate}
    \item Classical solver is used to produce dataset of pairs $(\text{PDE data}, \text{solution})$.
    \item Neural network $\mathcal{N}_1$ is trained on the dataset to predict solutions from PDE data.
    \item Predictions of neural network $\mathcal{N}_1$ are analyzed with the help of the upper bound, i.e., for each prediction certificate $w_{U}$ is obtained by direct optimization of the upper bound. From this data, a second dataset of triples $(\text{PDE data}, \text{solution}, \text{certificate }w_{U})$ is formed.
    \item Second neural network $\mathcal{N}_2$ is trained to predict certificate $w_{U}$ from PDE data and approximate solution produced by $\mathcal{N}_1$.
\end{enumerate}
At the inference stage, $\mathcal{N}_1$ predicts approximate solution $\widetilde{u}$ given problem data, $\mathcal{N}_2$ predicts certificate $w_{U}$ given approximate solution and problem data and the upper bound is used to compute upper bound on the error.

\subsection{Concrete error bounds}
\label{Section: Functional a posteriori error estimate as a novel loss function. Subsection: Concrete error bounds for elliptic PDEs.}

To illustrate abstract approach described in \cref{Section: Functional a posteriori error estimate as a novel loss function. SUbsection: High-level description of the error estimate.} we consider elliptic equation in the domain $\Gamma \subseteq [0, 1]^{D}$:
\begin{equation}
    \label{eq:general_elliptic_equation}
    \begin{split}
        -&\sum_{ij=1}^{D}\frac{\partial}{\partial x_i}\left(a_{ij}(x) \frac{\partial u}{\partial x_j} \right) + b^{2}(x) u(x) = f(x),\\
        &\left.u\right|_{\partial\Gamma} = 0,\,a_{ij}(x) \geq c > 0.
    \end{split}
\end{equation}
For later use we define several norms
\begin{equation}
    \label{eq:norms}
    \begin{split}
        &\left\|v\right\|_{2}^2 = \int_{\Gamma} dx\,v^2,\,\left\|w\right\|_{a^{-1}}^2 = \int_{\Gamma} dx\left(\sum_{i,j}\left(a^{-1}\right)_{ij} w_{i} w_{j}\right),\\
        &\vertiii{u}^2 = \left\|u\right\|^{2}_{a}  + \left\|b u\right\|_2^{2}.
    \end{split}
\end{equation}
Using the weak form, Cauchy-Schwarz, and Friedrichs inequalities, it is possible to show (see \cite{mali2013accuracy} (Chapter 3) and also \cref{Section: Derivation of Astral loss for elliptic equation}) that for \cref{eq:general_elliptic_equation} the energy norm of the deviation of approximate solution $\widetilde{u}$ from the exact one $u$ is bounded from above:
\begin{equation}
    \label{eq:upper_bound_nD}
    \begin{split}
        &\vertiii{\widetilde{u} - u_{\text{exact}}}^2 \leq \mathcal{L}_{\text{Astral}}[\widetilde{u}, y, \mathcal{D}, \beta] = \\
        &\int_{\Gamma} dx\frac{C^2(1+\beta)}{C^2 b(x)^2(1+\beta) + 1} \mathcal{R}(\widetilde{u}, y)^2 + c\left\|a\nabla \widetilde{u} - y\right\|^2_{a^{-1}},\\
        &\mathcal{R}(\widetilde{u}, y) = f(x) - b(x)^2\widetilde{u} + \sum_{i}\frac{\partial y_i(x)}{\partial x_i},\,c=\frac{1+\beta}{\beta},\\
        &C = 1 \big/\left(\inf_{x} \sqrt{\lambda_{\min}(a_{ij}(x))}\pi D\right),
    \end{split}
\end{equation}
where $\beta$ is positive number, $y$ is a vector field. Note that upper bound is the Astral loss defined in \cref{Section: Functional a posteriori error estimate as a novel loss function. Subsection: How to incorporate error estimates into deep-learning pipelines.} where $y$ and $\beta$ in \cref{eq:upper_bound_nD} correspond to $w_{U}$ from the definition \cref{eq:ASTRAL}.

\cref{eq:upper_bound_nD} possess properties described in \cref{Section: Functional a posteriori error estimate as a novel loss function}. First, it is easy to see that the bound is computable from the problem data (for considered PDE \cref{eq:general_elliptic_equation} the data is $b(x)$, $a_{ij}(x)$, $f(x)$), approximate solution $\widetilde{u}(x)$ and supplementary parameters (vector field $y$ and $\beta$). Second, it is straightforward to check that the bound is tight, since it is saturated when $y(x) = \sum_{j=1}^{D}a_{ij}\frac{\partial}{\partial x_{j}} u_{\text{exact}}$ and $\beta \rightarrow \infty$. For the formal proofs we refer to \cite{mali2013accuracy} (Chapter 3).

Our second example of functional a posteriori upper bound is given for initial-boundary value problem
\begin{equation}
    \label{eq:convection-diffusion}
    \begin{split}
        &\frac{\partial u(x, t)}{\partial t} - \frac{\partial^2 u(x, t)}{\partial x^2} + a\frac{\partial u(x, t)}{\partial x} = f(x, t),\\
        &u(x, 0) = \phi(x),\,\left.u(x, t)\right|_{x \in \partial \Gamma} = 0,\,a=\text{const},
    \end{split}
\end{equation}
where $\Gamma = [0, 1]\times[0, T]$. As shown in \cite{repin2010posteriori}, the upper bound for this case reads
\begin{equation}
    \label{eq:convection-diffusion-Astral}
    \begin{split}
        &\vertiii{e}_{\text{c.d.}} \equiv \int dx dt\left(\frac{\partial e}{\partial x}\right)^2 + \frac{1}{2} \int dx\,\left.e^2\right|_{t=T} \leq \mathcal{L}_{\text{Astral}}=\\
        &\int dx dt\left(\left(y - \frac{\partial v}{\partial x}\right)^2 + \frac{1}{\pi} \left(f - \frac{\partial v}{\partial t} - a \frac{\partial v}{\partial x} + \frac{\partial y}{\partial x}\right)^2\right),
    \end{split}
\end{equation}
where $e(x, t) = u(x, t) - v(x, t)$ is the error, $v(x, t)$ and $u(x, t)$ are approximate and exact solutions, $y(x, t)$ is a free field corresponding to $w_{U}$ from the definition \cref{eq:ASTRAL}.

We have seen that for elliptic and convection-diffusion equations, the bounds are available. It is also positive to derive similar bounds for other practically-relevant PDEs. In particular, bounds are known for Maxwell equations \cite{a223344}, reaction-diffusion problems \cite{MR2253056}, elastoplasticity problems \cite{FUNCTIONAL2009REPIN}, PDEs for the flow of viscose fluid \cite{Srepin2002posteriori}, and non-linear elliptic problems \cite{A_posteriorirepin1999posteriori}.

Besides, there are a few general techniques to obtain functional error estimates. Most notable is the use of dual theory for variational problems with convex functionals  \cite{repin2000posteriori}, \cite{repin2000posteriori1}. More details on functional error estimates and systematic ways to derive them are available in monographs \cite{repin2008posteriori}, \cite{neittaanmaki2004reliable}, \cite{mali2013accuracy}.

\subsection{Relation between energy and \texorpdfstring{$L_2$}{L2} norms}
Since the $L_2$ norm is the most widely used one, it is instructive to establish how it is related to the energy norm used in the previous section. For example, for elliptic equations, we can find
\begin{equation}
    \begin{split}
        &\vertiii{e}^2 = \int_{\Gamma} dx \sum_{ij} \frac{\partial e}{\partial x_{i}}\frac{\partial e}{\partial x_{j}}a_{ij} + \left\|be\right\|_2^2 \geq\\
        &\lambda_{\min}^2 \left\|e\right\|_2^2 + \left\|be\right\|_2^2 \geq \lambda_{\min}^2 ||e||^2 + \inf_{x} b(x)^2\left\|e\right\|_2^2,
    \end{split}
\end{equation}
where $\lambda_{\min}$ is a minimal eigenvalue of the elliptic problem $-\frac{\partial}{\partial x_{i}} a_{ij} \frac{\partial}{\partial x_{j}}u(x) = \lambda u(x)$ defined on the same domain and with the same boundary conditions as the original elliptic problem. From the expression above, we obtain the bound $\left\|e\right\|_2^2 \leq \frac{1}{\lambda_{\min}^2 + \inf_{x} \left(b(x)\right)^2} \vertiii{e}.$ With the same reasoning, we can obtain a lower bound and find that $L_2$ and energy norms are equivalent (for suitably defined space of functions). The upper bound alone is sufficient to claim that whenever the error energy norm is sufficiently small, the $L_2$ norm of error is small too.

\section{PiNN training with Astral, residual and variational losses}
\label{Section: Comparison of losses for physics-informed neural networks}

\begin{table*}
  \caption{Quality of the upper bound $\left(\sqrt{\mathcal{L}_{\text{Astral}}} - \vertiii{e}\right) \big/\vertiii{e}$ and comparison of error in two different norms for elliptic \cref{eq:general_elliptic_equation} in the $L-$shaped domain and convection-diffusion \cref{eq:convection-diffusion} equations.}
  \label{table:elliptic and convection diffusion}
  \centering
  \begin{tabular}{crrrr}
\toprule
 & \multicolumn{2}{c}{elliptic, L-shaped domain} & \multicolumn{2}{c}{convection-diffusion}\\
bound quality & \multicolumn{2}{c}{$0.84\pm0.60$} & \multicolumn{2}{c}{$0.63\pm0.18$} \\
\cmidrule(r){2-5}
loss & relative error &  energy norm & relative error &  energy norm \\
\midrule
Astral & $\left(7.3 \pm 9.8\right)10^{-3}$ & $\left(2.7 \pm 4.1\right)10^{-2}$ & $\left(8.3 \pm 8.9\right)10^{-3}$ & $\left(1.3 \pm 0.6\right)10^{-2}$ \\
residual & $\left(5.0 \pm 7.2\right) 10^{-3}$ & $\left(1.4 \pm 1.2\right)10^{-2}$ & $\left(9.9 \pm 11\right)10^{-3}$ & $\left(7.8 \pm 2.4\right)10^{-3}$\\
\bottomrule
\end{tabular}
\end{table*}

\begin{figure*}
  \centering
  \includegraphics[scale=0.65]{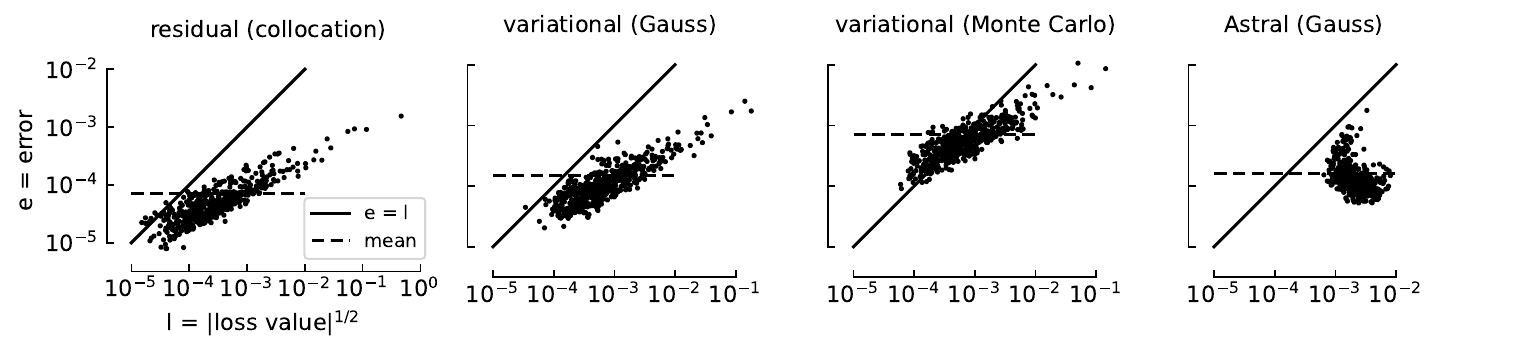}
  \caption{Comparison of PiNN solutions obtained with four different loss functions (from left to right): residual loss (\cref{eq:residual PiNN loss}, first equation) enforced on the set of points, variational loss (\cref{eq:residual PiNN loss}, second equation) computed with Gauss quadrature, variational loss computed with Monte Carlo, Astral loss \cref{eq:functional PiNN loss} computed with Gauss quadrature. On each plot $x$-axis represents square root of the absolute value of the loss and $y$-axis represents the value of error in the energy norm $\vertiii{e}$ (see \cref{eq:norms}), both computed at the end of the optimization process.}
  \label{fig:PiNN losses}
\end{figure*}

To practically evaluate Astral loss function, we compare it to other available losses for elliptic PDE. For simplicity we consider \cref{eq:general_elliptic_equation} with $a_{ij}(x) = a(x)\delta_{ij},\,b(x) = 0, D=2,$ where $\delta_{ij}$ are elements of identity matrix. For this particular case, upper bound \cref{eq:upper_bound_nD} simplifies (see \cref{Section: Simplified Astral loss for scalar case}), and we obtain a loss function
\begin{equation}
    \label{eq:functional PiNN loss}
    \begin{split}
    &\mathcal{L}_{\text{Astral}}[u, y, \mathcal{D}, \beta] = (1 + \beta)\\ &\int dx\left( \frac{\left(f(x) + \sum\limits_{i=1}^{2}\frac{\partial y_i}{\partial x_i}\right)^2}{4\pi^2\inf_{x}a(x)} + \frac{\sum\limits_{i=1}^{2}\left(a(x)\frac{\partial u}{\partial  x_i} - y_i\right)^2}{\beta a(x)}\right).
    \end{split}
\end{equation}
For the elliptic equation, one can also apply standard residual and variational loss functions\footnote{For consistency, here we write residual loss in a continuous form, but in our experiments, we used the standard collocation approach.}
\begin{equation}
    \label{eq:residual PiNN loss}
    \begin{split}
    &\mathcal{L}_{\text{res}}=\int dx \left(\sum_{i=1}^{2}\frac{\partial}{\partial x_{i}}\left( a\frac{\partial u}{\partial x_{i}}\right) + f\right)^2,\\
    &\mathcal{L}_{\text{var}} =\int dx \left(\frac{1}{2}a(x)\sum_{i=1}^{2}\left(\frac{\partial u}{\partial x_{i}}\right)^2 - fu\right),
    \end{split}
\end{equation}
widely used in physics-informed neural networks \cite{lagaris1998artificial}, \cite{raissi2019physics}, and Deep Ritz method \cite{yu2018deep}.

To compare three loss functions, we produce a dataset for elliptic equations with known exact solutions and randomly generated smooth $a(x)$ with widely varying magnitudes (see \cref{Appendix: Experiment 1: comparison of losses for physics-informed neural networks} for details on dataset and training). On this dataset, we train physics-informed networks of similar size under the same optimization setting for three different losses above. To approximate integral, we use Monte Carlo \cite{kalos2009monte} (Section 4) method or Gauss quadratures \cite{tyrtyshnikov1997brief} (Lecture 16), and in case of residual loss, we used standard collocation formulation. In the network design and training, we closely follow the best practices for PiNN training \cite{wang2023expert}. A complete description of the training, network architectures, and dataset generation are available in \cref{Appendix: Experiment 1: comparison of losses for physics-informed neural networks.}. To measure the quality of predicted solution we use error in the energy norm $\vertiii{e}$ (see \cref{eq:norms}).

In addition we perform two more experiments. Namely, training with residual and Astral loss \cref{eq:convection-diffusion-Astral} for convection-diffusion problem \cref{eq:convection-diffusion}, as well as for the elliptic problem \cref{eq:general_elliptic_equation} in $L$-shaped domain, i.e., $\Gamma = [0, 1]^{2}\backslash[0.5, 1]^2$. Details on these experiments are given in \cref{Appendix: Experiment 1: convection-diffusion problem} and
\cref{Appendix: Experiment 1: elliptic equation L-shaped domain} respectively. The $L-$shaped domain is included as a classical benchmark used in literature on a posteriori error estimate.

\begin{figure*}[!t]
    \centering
    \includegraphics[scale=0.3]{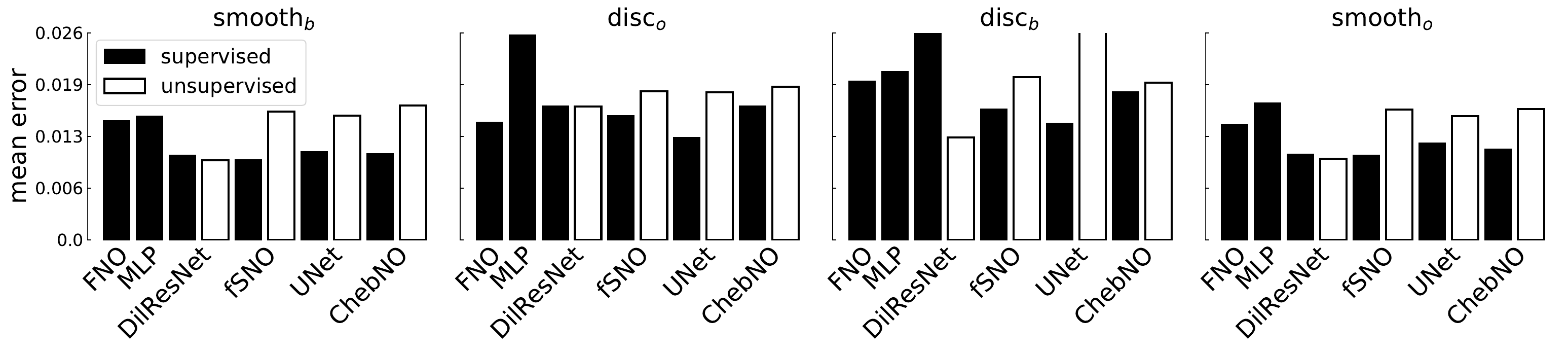}
    \caption{Comparison of mean errors in energy norm obtained with unsupervised (\cref{fig:unsupervised UQ training}) and supervised (\cref{fig:supervised UQ training}) training. Supervised training leads to slightly smaller error.}
    \label{fig:supervised vs unsupervised error}
\end{figure*}

\subsection{Discussion of training results}

Results for losses comparison appear in \cref{fig:PiNN losses}. We can make several important observations:
\begin{enumerate}
    \item For all losses but $\mathcal{L}_{\text{Astral}}$, there are samples with errors larger and smaller than loss (there are samples from both sides on the line with slope equals one). This means only $\mathcal{L}_{\text{Astral}}$ provides the upper bound for the error.
    \item Astral loss has a much smaller spread for the value of the loss and for the error, i.e., the optimization process with this loss is more robust. This can be explained by the fact, that Astral loss is equivariant to rescaling. Suppose we multiply elliptic equation \eqref{eq:general_elliptic_equation} by $s$. It is easy to see that for Astral loss \cref{eq:functional PiNN loss} one obtains overall scale $s^2$, i.e., the same scale as for squared error, simply by rescaling neural network output, i.e., $y_{i} \longrightarrow s y_{i}$. This is automatically done during training, so Astral loss converges more uniformly.
    \item On average, we see that residual loss shows faster convergence. However, when the number of iterations is increased, we reach roughly the same errors with all losses.
\end{enumerate}

\begin{figure}[!t]
    \centering
    \includegraphics[scale=0.3]{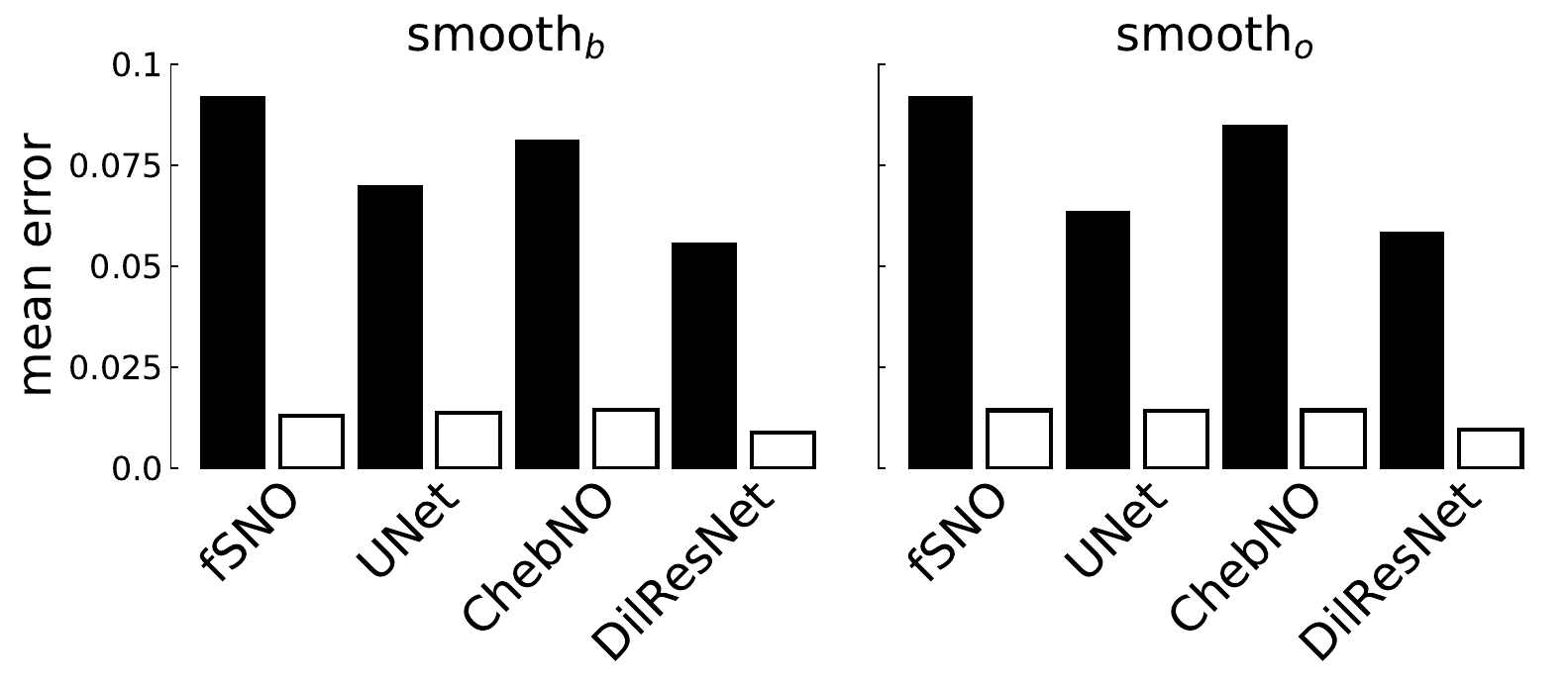}
    \caption{Comparison of errors for training with residual loss \cref{eq:residual PiNN loss} (black) used in PINO \cite{li2021physics} and Astral loss \cref{eq:functional PiNN loss} (white) for elliptic equations with smooth coefficient. The results for Astral loss are better by large margin.}
    \label{fig:PINO vs unsupervised upper bound}
\end{figure}

The results for convection-diffusion equation and elliptic equation in the $L-$
shaped domain appear in \cref{table:elliptic and convection diffusion}. We can highlight the following:
\begin{enumerate}
    \item Statistically, astral and residual losses typically provide roughly the same performance. On average, residual loss leads to slightly better relative error for elliptic problem and slightly worse for convection-diffusion equation.
    \item For both equations Astral loss provide high-quality upper bound. This is especially interesting for $L-$shaped domain since solution has corner singularity.
\end{enumerate}

\section{Application to parametric PDEs}
\label{Section: Application to parametric PDEs}

In the case of parametric PDE, our experiments were guided by the following research questions: (i) What is the difference in accuracy between supervised (\cref{fig:supervised UQ training}) and unsupervised (\cref{fig:unsupervised UQ training}) settings, which one leads to tighter bounds on error? (ii) How sensitive is training to the parameters of PDE \cref{eq:general_elliptic_equation}? In particular, how well does the training work for discontinuous diffusion coefficients and zero source terms? (iii) In the unsupervised setting (\cref{fig:unsupervised UQ training}), how does the upper bound loss compare with classical residual losses? (iv) How do the results depend on the choice of architecture?

To answer these questions, we generated four datasets with different source terms and diffusion coefficients. Namely, the equations with their short names are: $\text{smooth}_b$ -- $a(x)$ is smooth, $b(x)$ is present, \eqref{eq:elliptic_1}; $\text{disc}_{o}$ -- $a(x)$ is discontinuous, $b(x)$ is zero, \eqref{eq:elliptic_3}; $\text{disc}_{b}$ -- $a(x)$ is discontinuous, $b(x)$ is present, \eqref{eq:elliptic_4}; $\text{smooth}_{o}$ -- $a(x)$ is smooth, $b(x)$ is zero, \eqref{eq:elliptic_6}.

For these datasets we trained several architectures with different losses, including $L_2$ loss for the supervised training scheme (\cref{fig:supervised UQ training}), upper bound loss \cref{eq:upper_bound_nD} for $D=2$ and the residual loss \cref{eq:residual PiNN loss} for the unsupervised training scheme (\cref{fig:unsupervised UQ training}). The architectures used are: (i) FNO -- classical Fourier Neural Operator from \cite{li2020fourier}; (ii) fSNO -- Spectral Neural Operator on Gauss-Chebyshev grid. The construction mirrors FNO, but instead of FFT, a transformation based on Gauss quadratures is used \cite{fanaskov2022spectral}; (iii) ChebNO -- Spectral Neural Operator on Chebyshev grid. Again, the construction is the same as for FNO, but DCT-I is used instead of FFT \cite{fanaskov2022spectral}; (iv) DilResNet -- Dilated Residual Network from \cite{yu2017dilated}, \cite{stachenfeld2021learned}; (v) UNet -- classical computer vision architecture introduced in \cite{ronneberger2015u}; (vi) MLP -- vanilla feedforward neural network. So, we cover modern neural operators (FNO, fSNO, ChebNO) and classical computer vision architectures widely used for PDE modeling (DilResNet, UNet). MLP appears as a weak baseline.

To measure the quality of the predicted solution, we, again, use error in the energy norm $\vertiii{e}$ (see \cref{eq:norms}), and to measure the quality of predicted certificates $y$, $\beta$, we use the ratio of upper bound (the square root of the loss function \cref{eq:functional PiNN loss}) to the error in the energy norm. The choice of metric for certificates is motivated by the fact that the only end of $y$ and $\beta$ is to provide an upper bound. Also, only relative value is of interest since even when the error is large, neural network can provide a good upper bound, if it predicts the magnitude of this error accurately.

More details on experiments appear in \cref{Appendix: Experiment 2: supervised training with build-in error estimate.} and \cref{Appendix: Experiment 3: unsupervised training with build-in error estimate.}. More experimental results can be found in \cref{Appendix: Experiment 1: comparison of losses for physics-informed neural networks.}

\begin{figure*}[!t]
    \centering
    \includegraphics[scale=0.3]{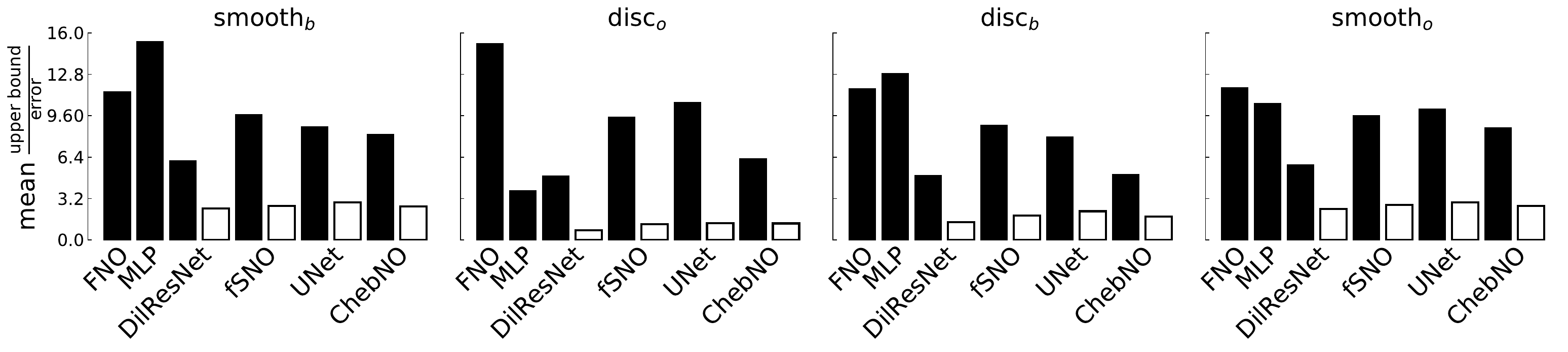}
    \caption{Comparison of relative upper bounds obtained with unsupervised (\cref{fig:unsupervised UQ training}) and supervised (\cref{fig:supervised UQ training}) training. Unsupervised training leads to much tighter error upper bound.}
    \label{fig:supervised vs unsupervised upper bound}
\end{figure*}
\subsection{Discussion of training results}
We can make the following important observations:
\begin{enumerate}
    \item We can see from \cref{fig:PINO vs unsupervised upper bound} that the Astral loss function produces much better error than the PINN loss \cite{li2021physics} for all architectures considered and all equations with smooth coefficients.
    \item \cref{fig:supervised vs unsupervised error} indicates that the supervised training scheme leads to slightly better errors for all architectures but DilResNet.
    \item From \cref{fig:supervised vs unsupervised upper bound}, we can conclude that unsupervised training leads to a much tighter error bound than the supervised one. The upper bound obtained with Astral loss is almost the same as the one obtained with numerically expensive direct optimization of the upper bound. This shows one can produce high-quality upper bounds cheaply with a neural network.
    \item We do not observe any sensitivity to the smoothness of coefficients of the elliptic equation, the presence of the source term, and the architecture used. In general, DilResNet performs better than other architectures, and neural operators seem to show second-best results, but the differences are not prominent.
\end{enumerate}

\section{Related research}
\label{Section: Related research}
We can classify the related research articles into several categories.

The first one consists of generalization error analysis for physics-informed neural networks, e.g., \cite{de2022error}, \cite{mishra2022estimates}, \cite{de2022error_1}, \cite{gonon2022uniform}, \cite{jiao2021error}. In such contributions, authors show that it is possible to reach a given approximation accuracy with a neural network of a particular size. The bound on error involves unknown quantities, such as the norm of the exact solution. This line of work provides scaling arguments and establishes how expressive a particular architecture can be, but the results are unusable for actual error estimates.

The second line of work is on adaptive loss estimation for physics-informed neural networks. The examples include \cite{wu2023comprehensive}, \cite{zubov2021neuralpde}. In a standard physics-informed neural network, the set of collocation points is either randomly selected \cite{yu2018deep} or fixed from the start \cite{raissi2019physics}. Adaptive loss estimation aims to select points such that the residual is uniformly small. Typically, this tends to improve the accuracy and stability of training. The problem is it does not allow for error control since the relation between the residual and the error is not straightforward.

The next set of articles is related to the generalization bound for neural operators, e.g., \cite{kovachki2021universal}, \cite{lanthaler2022error}. In analogy with physics-informed neural networks, the results usually established are bounds on network size sufficient to solve a given class of parametric PDEs. For the same reasons, with those results, it is not possible to estimate the actual error in a practical setting.

Besides that, there are several contributions directly related to a posteriori error estimation \cite{guo2022energy}, \cite{filici2010error}, \cite{hillebrecht2022certified}, \cite{berrone2022solving}, \cite{minakowski2023priori}, \cite{roth2022neural}, \cite{cai2020deep}. In all these contributions, authors specialize in a particular classical error bound to deep learning PDE/ODE solvers. Namely, \cite{filici2010error} adopts a well-known error estimation for ODEs, based on the construction of related problems with exactly known solution \cite{zadunaisky1976estimation}. Similarly, \cite{hillebrecht2022certified} uses well-known exponential bound on error that involves residual and Lipschitz constant \cite{hairer1993solving} and applies a neural network to perform the residual interpolation. Similarly, contributions \cite{guo2022energy}, \cite{berrone2022solving} and \cite{cai2020deep} are based on FEM posterior error estimates, and \cite{minakowski2023priori}, \cite{roth2022neural} are on dual weighted residual estimator \cite{becker2001optimal}.

\section{Conclusion}
\label{Section: Conclustion}
Following \cite{Muzalevsky2021aposteriori}, we outlined how a posteriori error estimate of functional type can be used to train physics-informed neural networks, and neural operators in supervised and unsupervised settings. As we argue, error majorants provide a novel systematic way to construct loss functions for physics-informed training. Our results indicate that for elliptic equations, the proposed loss leads to more stable training than classical residual and variational losses. Besides that, it allows for direct error control. For the parametric PDE, functional a posteriori error estimate allows for the first time to obtain practical error upper bound. That way, one can estimate the quality of the solution obtained by a black-box neural network. Predicting error majorants with neural networks looks especially appealing since one simultaneously retains guarantees (thanks to the superb properties of functional error estimate) and substantially decreases the computational load by removing the need to solve a complex optimization problem to obtain a tighter upper bound.

\section{Reproducibility}
\label{Section: Reproducibility}
We use JAX \cite{jax2018github}, Equinox \cite{kidger2021equinox} and Optax \cite{deepmind2020jax}.

GitHub repository with notebooks, architectures and training scripts is \href{https://github.com/VLSF/UQNO}{https://github.com/VLSF/UQNO}

\section{Acknowledgements}
The authors would like to thank Anna Rudenko, who suggested the name Astral (neur\textbf{A}l a po\textbf{ST}erio\textbf{R}i function\textbf{A}l \textbf{L}oss) for the loss function \eqref{eq:ASTRAL}.

\bibliography{refs.bib}
\bibliographystyle{icml2024}

\newpage
\appendix
\onecolumn

\section{Derivation of Astral loss for elliptic equation}
\label{Section: Derivation of Astral loss for elliptic equation}
The discussion in this section follows Chapter 3 of \cite{mali2013accuracy}.

We consider the boundary value problem
\begin{equation}
\begin{aligned}
- \text{div} ~a \nabla u(x) &+ b(x)^2 u(x) = f(x), \quad x \subseteq \Gamma = [0, ~1]^D, \\
&u\big|_{\partial\Gamma} = 0.
\label{eq:bv}
\end{aligned}
\end{equation}
where $a$ is a symmetric matrix satisfying the condition $ax \cdot x \ge c \vert x\vert^2, \quad \forall x \in \mathbb{R}^D$  and  $b^2$ is a nonnegative function of $x$. The general solution $u(x)$ of \eqref{eq:bv} is defined by the integral identity
\begin{equation}
\int_{\Gamma}dx \Big[a \nabla u(x) \cdot \nabla w(x) + b^2(x) u(x) w(x)\Big] = \int_{\Gamma}dx ~f(x)w(x),
\label{eq:weak_form}
\end{equation}
that holds for every $w(x)$ from Sobolev space of square summable functions with square summable derivatives up to the order 1,~$H^1([0, ~1]^D)$, with Dirichlet boundary conditions. Let $\widetilde{u}(x)$ approximate the exact solution $u(x)$  of the problem \eqref{eq:bv}. By \eqref{eq:weak_form} we deduce the relation
$$
\int_{\Gamma}dx \Big[a \nabla(u - \widetilde{u}) \cdot \nabla w + b^2 (u - \widetilde{u}) w\Big] = \int_{\Gamma}dx \Big[fw - b^2 \widetilde{u}w - \nabla \widetilde{u} \cdot \nabla w\Big].
$$

Since $w$ vanishes at the boundary:
$$\int_{\Gamma} dx (\nabla \cdot (yw)) = yw\Big|_{\partial \Gamma} = 0,$$
where $y(x) \in \mathbb{R}^D$ is arbitrary vector. We can rewrite previous relation as follows:
\begin{equation}
\int_{\Gamma}dx \Big[a \nabla(u - \widetilde{u}) \cdot \nabla w + b^2 (u - \widetilde{u}) w\Big] = \int_{\Gamma}dx \Big[\mathcal{R}(\widetilde{u}, ~y) w + \big(y - a \nabla \widetilde{u}\big) \cdot \nabla w\Big],
\label{eq:mod_wf}
\end{equation}
where  $\mathcal{R}(\widetilde{u}, y) = f - b^2 \widetilde{u} + \text{div}~y$. Let us represent the first integral on the right-hand side of \eqref{eq:mod_wf} as follows
$$\int_{\Gamma}dx ~\mathcal{R}(\widetilde{u}, ~y)w = \int_{\Gamma} dx ~\alpha \mathcal{R}(\widetilde{u}, ~y)w + \int_{\Gamma} dx ~(1 - \alpha)\mathcal{R}(\widetilde{u}, ~y)w,$$
where $\alpha \in L_{[0,~1]}^{\infty}\big(\Gamma\big) = \big\{\alpha \in L^{\infty}\big(\Gamma\big)~\big|~ 0 \le \alpha(x) \le 1\big\}$ is a weight function. It is easy to see that
$$
\int_{\Gamma}dx \big(y - a \nabla \widetilde{u}\big) \cdot \nabla w \le \big\Vert y - a \nabla \widetilde{u}\big\Vert_{a^{-1}} \Vert \nabla w\Vert_a,
$$
$$
\int_{\Gamma}dx ~\mathcal{R}(\widetilde{u}, ~y)w \le C \Big\Vert \mathcal{R}(\widetilde{u}, ~y)\Big\Vert_2 \Vert \nabla w\Vert_a,
$$
where $C$ is a constant in the inequality
$$
\Vert w\Vert_2 \le C \Vert \nabla w\Vert_a, \quad \forall w \in H^{1}\big([0,~1]^D\big).
$$

Then, we have
\begin{equation}
\begin{aligned}
\Big\vert \int_{\Gamma} dx~\mathcal{R}(\widetilde{u},~y)w\Big\vert &\le \Big\Vert \frac{\alpha}{b} \mathcal{R}(\widetilde{u},~y)\Big\Vert_2 \Vert bw\Vert_2 + C \Big\Vert (1 - \alpha) \mathcal{R}(\widetilde{u},~y)\Big\Vert_2 \Vert \nabla w\Vert_a, \\
& \Vert w \Vert_a^2 = \int_{\Gamma}dx~ aw \cdot w, \quad \Vert w \Vert_{a^{-1}}^2 = \int_{\Gamma}dx~ a^{-1}w \cdot w.
\end{aligned}
\end{equation}

By  setting $w = u - \widetilde{u}$  we arrive at the estimate
$$
\vertiii{\widetilde{u} - u}^2  \le \Big(C \Big\Vert \big(1 - \alpha \big) \mathcal{R}(\widetilde{u}, ~y)\Big\Vert_2 + \big\Vert a\nabla \widetilde{u} - y \big\Vert_{a^{-1}}\Big)^2 + \Big\Vert \frac{\alpha}{b}\mathcal{R}(\widetilde{u}, ~y)\Big\Vert^2_2.
$$

In the sake of simplicity, we use the algebraic Young’s inequality. For $a,~b \in \mathbb{R}$ and for any $\beta$ positive number:
$$
2ab \le \beta a^2 + \frac{1}{\beta}~b^2.
$$
We find that
$$
\big\vert a + b \big\vert^2 \le \big(1 + \beta\big) \big\vert a \big\vert^2 + \frac{1 + \beta}{\beta}\big\vert b \big\vert^2.
$$

Thus, we have
\begin{equation}
\vertiii{\widetilde{u} - u}^2 \le (1 + \beta) C^2 \Big\Vert\big( 1 - \alpha \big) \mathcal{R}(\widetilde{u}, ~y)\Big\Vert^2_2 + \frac{1+\beta}{\beta} \big\Vert a\nabla \widetilde{u} - y\big\Vert_{a^{-1}} + \Big\Vert \frac{\alpha}{b}\mathcal{R}(\widetilde{u}, ~y)\Big\Vert^2_2,
\label{eq:ub}
\end{equation}
where $\beta$  is an arbitrary positive number.

Minimization of the right-hand side of \eqref{eq:ub} with respect to $\alpha$ is reduced to the following auxiliary variational problem: find $\widehat{\alpha} \in L_{[0,~1]}^{\infty}\big(\Gamma\big)$ such that
$$
\Upsilon \big(\widehat{\alpha}\big) = \inf_{\alpha \in L_{[0,~1]}{\infty}(\Gamma)} \int_{\Gamma}dx~\Big(\alpha^2 P(x) + \big(1 - \alpha\big)^2 Q(x)\Big),
$$
and $P$ and $Q$ are nonnegative integrable functions, which do not vanish simultaneously. It is easy to find that for almost all $x$
$$
\widehat{\alpha}(x) = \frac{Q}{P + Q} \in [0,~1], \quad \Upsilon \big(\widehat{\alpha}\big) = \frac{PQ}{P + Q}.
$$
In our case, $P = b^{-2} \mathcal{R}\big(\widetilde{u},~y\big)$ and $Q = C^2 \big(1 + \beta\big) \mathcal{R}\big(\widetilde{u},~y\big).$

Therefore, we obtain
\begin{equation}
\begin{aligned}
\vertiii{\widetilde{u} - u}^2 &\le \int_{\Gamma} dx\frac{C^2(1+\beta)}{C^2 b(x)^2(1+\beta) + 1} \mathcal{R}(\widetilde{u}, y)^2 + \frac{1+\beta}{\beta}\left\|a\nabla \widetilde{u} - y\right\|^2_{a^{-1}} = \mathcal{L}_{\text{Astral}}[\widetilde{u}, y, \mathcal{D}, \beta], \\
&\mathcal{R}(\widetilde{u}, y) = f(x) - b(x)^2\widetilde{u} + \text{div}~y,\,C = 1 \big/\left(\inf_{x} \sqrt{\lambda_{\min}(a)}\pi D\right).
\end{aligned}
\end{equation}

\section{Simplified Astral loss for scalar diffusion coefficient}
\label{Section: Simplified Astral loss for scalar case}
In this section, we illustrate how to derive \cref{eq:functional PiNN loss} from \cref{eq:upper_bound_nD}.

First, observe that $b=0$ so the residual weight simplifies to $C^2(1+\beta)$ and the residual itself becomes $f(x) - \sum_{i=1}^{2}\frac{\partial y_{i}}{\partial x_{i}}$ for $D=2$ problem considered. Next, the second term with flux condition simplifies because $a_{ij}(x) = \delta_{ij} a(x)$ where $\delta_{ij}$ are matrix elements of the identity matrix. More specifically, the inverse matrix $a^{-1}$ is simply a division by $a(x)$, so the term becomes $\frac{1}{a(x)}\sum_{i=1}^{2}\left(a(x)\frac{\partial \widetilde{u}(x)}{\partial x_{i}} - y_i\right)^2$. Finally, the constant $C$ simplifies to $1 \big / \left(2 \pi \inf_{x}\sqrt{a(x)}\right)$, since for diagonal matrix $a_{ij}(x)$ the minimal eigenvalue is $a(x)$.

With all that we obtain a simplified version of Astral loss \cref{eq:functional PiNN loss}.

\section{Experiment 1: comparison of losses for physics-informed neural networks}
\label{Appendix: Experiment 1: comparison of losses for physics-informed neural networks}

In all experiments, we use feedforward physics-informed neural networks with $3$ layers, $50$ features in each layer, and $\text{GELU}$ activation functions. Following suggestions in \cite{wang2023expert}, we transform input coordinates to obtain $50$ Fourier features. Each network was optimized for $50000$ epoch with Lion optimizer \cite{chen2023symbolic}. Starting learning rate is $10^{-4}$ with exponential decay $\times 0.5$ each $10000$ epoch. The output of the neural network was multiplied by $\sin(\pi x)\sin(\pi y)$ to enforce boundary conditions. When error majorant is used as loss functions, we predict $y_1$, $y_2$, and $u$ with three separate neural networks.

To generate random elliptic equation we draw random function from $f \sim \mathcal{N}\left(0, \left(I-\Delta\right)^{-2}\right)$ after that this function is shifter and rescaled such that $\min_{x} f(x) = 1, \max_{x} f(x) = 6$. After that we multiply this function on $1/s$ where $s$ is drawn from exponential distribution with mean $100$. This way we obtain positive functions of widely varying scale that we use as diffusion coefficient $a(x)$. For all elliptic equations we fix exact solution to be $u(x) = x_1(1-x_1)x_2(1-x_2)$ and find $f(x)$ from known $a(x)$ and $u(x)$.

\section{Experiment 2: convection-diffusion problem}
\label{Appendix: Experiment 1: convection-diffusion problem}

The architecture and training scheme are the same as in \cref{Appendix: Experiment 1: comparison of losses for physics-informed neural networks}, but we train for $10000$ epoch and multiply learning rate by $0.5$ each $3000$ epoch.

To test PiNN networks we sample exact solutions of convection-diffusion equation \cref{eq:convection-diffusion} given by
\begin{equation}
    \label{eq:exact c-d}
    u(x, t) = \text{Re} \left(\sum_{k=0}^{N} c_ke^{-(2\pi k)^2 t - 2\pi k a i t + 2 \pi k x i}\right) \sin(\pi x) = \widetilde{u}(x, t) \sin(\pi x).
\end{equation}
It is easy to see that
\begin{equation}
    \frac{\partial \widetilde{u}(x, t)}{\partial t} - \frac{\partial^2 \widetilde{u}(x, t)}{\partial x^2} + a\frac{\partial \widetilde{u}(x, t)}{\partial x} = 0,
\end{equation}
and that $u(x, t)$ fulfills boundary conditions. If we choose
\begin{equation}
    \label{eq:f c-d}
    f(x, t) = \widetilde{u}(x, t)\left(\pi^2 \sin(\pi x) + a\pi \cos(\pi x)\right) - 2 \pi \frac{\partial \widetilde{u}(x, t)}{\partial x} \cos(\pi x),
\end{equation}
than $u(x, t)$ solves original equation \cref{eq:convection-diffusion}.

So, to generate exact solutions we used \cref{eq:exact c-d} with $c_k$ sampled from $\left(x + i y\right)\big/ (1 + (\pi k / 5)^2)^2$ where $x, y \sim \mathcal{N}(0, 1)$ and $a$ sampled from $10^{-2}\mathcal{N}(0, 1)$, $N=150$. After that we compute source term using \cref{eq:f c-d}.

\section{Experiment 3: elliptic equation, L-shaped domain}
\label{Appendix: Experiment 1: elliptic equation L-shaped domain}

The architecture and training scheme are the same as in \cref{Appendix: Experiment 1: convection-diffusion problem}.

For this particular geometry exact solutions are available only for specific boundary conditions, so we use finite difference discretization to obtain solution on the fine grid and use this solution as a ground truth.

Parameters of the problem are sampled as follows
\begin{equation}
    a(x) = \left(10 g_1(x)\right)^2 + 1,\,f(x) = 10 g_2(x),\, b(x)=0,\, g_1(x), g_2(x)\sim \mathcal{N}\left(0,\left(I - 10^{-1}\Delta\right)^{-2}\right),
\end{equation}
where Laplace operator has periodic boundary conditions.

To enforce Dirichlet boundary conditions in all experiments we used mean value coordinates \cite{floater2003mean} as explained in \cite{sukumar2022exact}.

\section{Elliptic equations used for experiments with parametric PDEs}

All elliptic PDEs considered are based on general form, described in \eqref{eq:general_elliptic_equation}.

First, we define random trigonometric polynomials
\begin{equation}
    \label{eq:random_trig_2D}
    \mathcal{P}(N_1, N_2, \alpha) = \left\{f(x) = \mathcal{R}\left(\sum_{m=0}^{N_1}\sum_{n=0}^{N_2}\frac{c_{mn}\exp\left(2\pi i(mx_1 + nx_2)\right)}{(1+m+n)^\alpha}\right):\mathcal{R}(c), \mathcal{I}(c)\simeq \mathcal{N}(0, I)\right\}.
\end{equation}

For the first equation, we use Cholesky factorization to define matrix $a$ and random trigonometric polynomials for $b$ and $f$:
\begin{equation}
    \label{eq:elliptic_1}
    \begin{split}
    &a(x) = \begin{pmatrix}
        \alpha(x) & 0 \\
        \gamma(x) & \beta(x)
    \end{pmatrix}
    \begin{pmatrix}
        \alpha(x) & \gamma(x) \\
        0 & \beta(x)
    \end{pmatrix},\\
    &\,\alpha(x),\beta(x) \simeq 0.1\mathcal{P}(5, 5, 2) + 1;\,\gamma(x),\,b(x),\,f(x) \simeq \mathcal{P}(5, 5, 2).
    \end{split}
\end{equation}

The next equation has a discontinuous scalar diffusion coefficient:
\begin{equation}
    \label{eq:elliptic_3}
    \begin{split}
    a(x) = \alpha(x)I,\,\alpha(x) = \begin{cases}
        10,\,p_1(x) \geq 0;\\
        1,\, ~~p_1(x) < 0,
    \end{cases}\,b(x) = 0,\,f(x) = 1,\,p_1(x) \simeq \mathcal{P}(5, 5, 2).
    \end{split}
\end{equation}
The analogous equation is dubbed ``Darcy flow'' in \cite{li2020fourier}.

Third equation is similar to \eqref{eq:elliptic_3} but with more diverse $b$ and $f$:
\begin{equation}
    \label{eq:elliptic_4}
    \begin{split}
    a(x) = \alpha(x)I,\,\alpha(x) = \begin{cases}
        10,\,p_1(x) \geq 0;\\
        1,\, ~~p_1(x) < 0,
    \end{cases}\,b(x),\,f(x),\,p_1(x) \simeq \mathcal{P}(5, 5, 2).
    \end{split}
\end{equation}

Last equation is similar to \eqref{eq:elliptic_1} but with $b=0$:
\begin{equation}
    \label{eq:elliptic_6}
    \begin{split}
    &a(x) = \begin{pmatrix}
        \alpha(x) & 0 \\
        \gamma(x) & \beta(x)
    \end{pmatrix}
    \begin{pmatrix}
        \alpha(x) & \gamma(x) \\
        0 & \beta(x)
    \end{pmatrix},\\
    &\,\alpha(x),\beta(x) \simeq 0.1\mathcal{P}(5, 5, 2) + 1;\,\gamma(x), f(x) \simeq \mathcal{P}(5, 5, 2);\,b(x) = 0.
    \end{split}
\end{equation}

\section{Experiment 4: supervised training with build-in error estimate}
\label{Appendix: Experiment 2: supervised training with build-in error estimate.}

We performed $6$ training runs with $N_{\text{train}}=200, 400, 600, 800, 1000, 1200$ for all architectures listed in \cref{Section: Application to parametric PDEs} and all equations. The number of grid points is $2^5 + 1$ along each dimension. To estimate the exact solution, we solve the same problem with a higher resolution grid, using $2^7 + 1$ points. The averaged results are reported in Tables~\ref{table:eq_2D},\,\ref{table:model_2D},\,\ref{table:N_train_2D},\,\ref{table:quality_2D}.

In addition, for the first equation \eqref{eq:elliptic_1}, we gather data on the same $6$ training runs with varying resolution $2^J + 1$, $J=5, 6, 7$. For this experiment, averaged results appear in Table~\ref{table:J_2D}.

In all optimization runs, we train for $1000$ epochs with Adam optimizer having the learning rate $10^{-3}$ multiplied by $0.5$ each $200$ epochs (exponential decay) and weight decay $10^{-2}$.

Details on architectures used are as follows.
\begin{enumerate}
    \item Construction of FNO closely follows the one given in \cite{li2020fourier}. We use $24$ featrures in processor, $4$ layers and the number of modes is $\left\lceil N_{\text{spatial}} / 4\right\rceil$.
    \item Construction of fSNO is similar to the one of FNO (i.e., encoder-processor-decoder and integral kernel in place of the linear layer), but following \cite{fanaskov2022spectral} we replace Fourier basis with orthogonal polynomials in the integral kernel. For this particular architecture, the construction of the integral kernel is as follows. We use the Gauss-Chebyshev grid, and compute projection on the basis of polynomials using Gauss quadratures \cite{golub1969calculation}. After the projection on the space of polynomials, we apply three convolutions with kernel size $3$. Finally, to return to the physical space, we compute sum $\sum_{n} c_n p_n(x)$ on the Gauss-Chebyshev grid, where $c_n$ are coefficients obtained after convolutions. In this case, the number of features in the processor is $34$, the number of layers is $4$ and the number of modes is $\left\lceil N_{\text{spatial}} / 4\right\rceil$
    \item ChebNO is a spectral neural operator \cite{fanaskov2022spectral} defined on Chebyshev grid. Construction of integral kernel is similar to fSNO but DCT is used to find projection on the polynomial space in place of Gauss quadratures. In the current architecture, the processor has $32$ features, the number of layers is $4$, and $16$ modes are used for all grids.
    \item DilResNet that we use closely follows architecture described in \cite{stachenfeld2021learned}. Namely, we use a processor with $24$ features that has $4$ layers. Each layer consists of convolutions with strides $[1, 2, 4, 8, 4, 2, 1]$, kernel size $3$, and skip connection.
    \item For UNet \cite{ronneberger2015u} we start with $10$ features and double the number of features with each downsampling that decreases the number of grid points by the factor of $2$ in each dimension. On each grid we use $2$ convolutions (kernel size $3$) and max pooling, transposed convolution are used for upsampling, and $3$ convolutions (kernel size $3$) appears on each grid after upsampling. In total, we have $4$ grids.
    \item MLP that we use consists of linear layers that process each dimension (including the feature dimension) separately. That way, the linear operator is defined by three matrices in $D=2$. MLP uses $64$ in the processor and has $4$ layers.
\end{enumerate}
For all networks, we use $\text{ReLU}$ nonlinearity.
Typical number of parameters for each network is given in the table below.

\begin{center}
\begin{tabular}{lcccccc}
\toprule
& FNO & fSNO & ChebNO & DilResNet & UNet & MLP \\
\# parameters& $668\times 10^{3}$ & $130\times 10^{3}$ & $115\times 10^{3}$ & $147\times 10^{3}$ & $248\times 10^3$ & $24\times 10^3$\\
\bottomrule
\end{tabular}
\end{center}

\section{Experiment 5: unsupervised training with build-in error estimate}
\label{Appendix: Experiment 3: unsupervised training with build-in error estimate.}

The unsupervised case is addressed by analyzing the $D = 2$ of \eqref{eq:general_elliptic_equation} and optimizing the upper bound \eqref{eq:upper_bound_nD} to produce an approximate solution $\widetilde{u}$ and an upper bound certificate $y$. The loss that was used to train the neural network $\mathcal{N}$ is:
\begin{equation}
    \label{eq:upper_bound_loss}
    \mathcal{L}[\widetilde{u}, y, \mathcal{D}, \beta, \lambda] = \sqrt{\mathcal{L}_{\text{Astral}}[\widetilde{u}, y, \mathcal{D}, \beta]} + \lambda \sqrt{\widetilde{u}^2_{\partial \Gamma}},
\end{equation}
\[
    \mathcal{L}[\widetilde{u}, y, \mathcal{D}, \beta, \lambda] \rightarrow \min_{\widetilde{u},~y},
\]
where $\widetilde{u},~y$ are the output of the neural network, $a, b, f, C$ are defined from problem data $\mathcal{D}$ for each equation, and $\lambda$  and $\beta$ are hyperparameters of the loss function \eqref{eq:upper_bound_loss}. Including the part of the loss  with boundaries is necessary to compensate for the lack of information about the exact solution on $x \in \partial\Gamma$.

In this paper, a comparison was made with the state-of-the-art unsupervised model proposed by \cite{li2021physics}. The PINO loss is defined as
\begin{equation}
    \label{eq:pino_loss}
    \mathcal{L}_{\text{pino}}[\widetilde{u}, u_{\text{exact}}, \mathcal{D}, \alpha, \gamma] =  \mathcal{L}_{\text{data}}[\widetilde{u}, u_{\text{exact}}] + \alpha \mathcal{L}_{\text{residual}}[\widetilde{u}, \mathcal{D}] + \gamma \sqrt{\widetilde{u}^2_{\partial \Gamma}},
\end{equation}
where
\begin{equation*}
    \begin{split}
        &\mathcal{L}_{\text{residual}}[\widetilde{u}, \mathcal{D}]= \sqrt{\int_{\Gamma}dx\Bigg(\sum_{i=1}^{2}\frac{\partial}{\partial x_{i}}\left( a(x)\frac{\partial \widetilde{u}}{\partial x_{i}}\right) + f(x) - b(x)\widetilde{u}\Bigg)^2}, \\
        &\mathcal{L}_{\text{data}}[\widetilde{u}, u_{\text{exact}}] = \sqrt{\int_{\Gamma}dx \Big(\widetilde{u} - u_{\text{exact}}\Big)^2}.
    \end{split}
\end{equation*}
The PINO loss consists of the physics loss in the interior and the data loss on the boundary condition, with hyperparameters $\alpha, \gamma > 0$. It is important to note that the proposed loss function \eqref{eq:upper_bound_loss} does not contain an exact solution $u_{\text{exact}}$ obtained using traditional solvers, unlike the PINO loss function \eqref{eq:pino_loss}.

We trained for $500$ epochs with Adam optimizer that has the learning rate $2 \cdot 10^{-3}$ multiplied by $0.5$ for each $50$ epoch and weight decay $10^{-2}$. Hyperparameters $\lambda$ and $\beta$ were chosen to be $1$.

The experiments are carried out according to the scheme in \cref{fig:unsupervised UQ training}. We performed $9$ training runs with $N_{\text{train}}=200, 400, 600, 800, 1000, 1200, 1400, 1600, 1800$ for all architectures and all equations. The number of points in the grid is $2^5 + 1$ along each dimension. To estimate the exact solution, we solve the same problem with a higher resolution grid, using $2^7 + 1$ points. The results using our loss function \eqref{eq:upper_bound_loss} are presented in Tables \ref{table:unsupervised_learning_energy_norm_for_all_models_1800},\ref{table:unsupervised_learning_fSNO_for_different_2D_datasets},\ref{table:unsupervised_learning_ChebNO_for_different_2D_datasets},\ref{table:unsupervised_learning_UNet_for_different_2D_datasets},\ref{table:unsupervised_learning_DilResNet_for_different_2D_datasets}. The results for PINO loss proposed in \cite{li2021physics} are presented in Table \ref{table:pino_loss_for_all_models_for_all_models_1800}.

The architectures used for unsupervised learning are identical to those used for supervised learning. For all networks, we use ReLU nonlinearity.

In Tables \ref{table:unsupervised_learning_fSNO_2D_for_different_lambda}, \ref{table:unsupervised_learning_ChebNO_2D_for_different_lambda}, \ref{table:unsupervised_learning_UNet_2D_for_different_lambda}, \ref{table:unsupervised_learning_DilResNet_2D_for_different_lambda}, one can find the results for train runs with different $\lambda$ for all architectures and datasets for 2D equations. To check dependencies from $\lambda$, we used  $N_{\text{train}}=1200$ for training, and $N_{\text{test}}=800$ for testing.


\section{Additional data for experiments 2 and 3}

\label{Appendix: Experiment 1: comparison of losses for physics-informed neural networks.}
\begin{table}[!th]
  \caption{Results for $N_{\text{train}} = 1800$, $N_{\text{test}} = 200$ for 2D equations using loss \eqref{eq:upper_bound_loss}.}
  \label{table:unsupervised_learning_energy_norm_for_all_models_1800}
  \centering
\begin{tabular}{cllllllll}
\toprule
{equation}  &                \multicolumn{2}{c}{fSNO} & \multicolumn{2}{c}{ChebNO} & \multicolumn{2}{c}{DilResNet} & \multicolumn{2}{c}{UNet} \\
{} &  $E_{\text{train}}$ &     $E_{\text{test}}$ &  $E_{\text{train}}$ &  $E_{\text{test}}$ &          $E_{\text{train}}$ &  $E_{\text{test}}$ &              $E_{\text{train}}$ &  $E_{\text{test}}$ \\

\midrule
\eqref{eq:elliptic_1}  &  0.012  &  0.013  &  0.013  &  0.015  &    0.008  &  0.009  &  0.012  &  0.014 \\
\eqref{eq:elliptic_3}  &  0.019  &  0.018  &  0.019  &  0.018  &    0.017  &  0.016  &  0.025  &  0.024 \\
\eqref{eq:elliptic_4}  &  0.014  &  0.017  &                      0.015  &  0.018  &  0.011  &  0.012  &                 0.070  &  0.078 \\
\eqref{eq:elliptic_6}   &   0.013  &  0.014  &                      0.013  &  0.014  &  0.009  &  0.010  &                  0.013  &  0.014 \\
\bottomrule
\end{tabular}
\end{table}

\begin{table}[!th]
  \caption{Results for $N_{\text{train}} = 1800$, $N_{\text{test}} = 200$ for 2D equations using PINO loss \eqref{eq:pino_loss}.}
  \label{table:pino_loss_for_all_models_for_all_models_1800}
  \centering
  \begin{tabular}{cllllllll}
\toprule
{equation}  &                \multicolumn{2}{c}{fSNO} & \multicolumn{2}{c}{ChebNO} & \multicolumn{2}{c}{DilResNet} & \multicolumn{2}{c}{UNet} \\
{} &  $E_{\text{train}}$ &     $E_{\text{test}}$ &  $E_{\text{train}}$ &  $E_{\text{test}}$ &          $E_{\text{train}}$ &  $E_{\text{test}}$ &              $E_{\text{train}}$ &  $E_{\text{test}}$ \\

\midrule
\eqref{eq:elliptic_1}   &                0.091 &               0.092 &                       0.081 &               0.081 &                       0.055 &               0.056 &                       0.065 &                0.070 \\
\eqref{eq:elliptic_6}   &                0.091 &               0.092 &                       0.082 &               0.085 &                       0.058 &               0.058 &                       0.061 &               0.064 \\
\bottomrule
\end{tabular}
\end{table}
\begin{table}[!th]
  \caption{Results for distinct $N_{\text{train}}$ averaged with respect to equations for 2D using loss \eqref{eq:upper_bound_loss}, $N_{\text{test}} = 200$.}
  \label{table:unsupervised_learning_for_different_train_size_for_all_models_averaged_wrt_equations}
  \centering
  \begin{tabular}{cllllllll}
\toprule
$N_{\text{train}}$  &                \multicolumn{2}{c}{fSNO} & \multicolumn{2}{c}{ChebNO} & \multicolumn{2}{c}{DilResNet} & \multicolumn{2}{c}{UNet} \\
{} &  $E_{\text{train}}$ &     $E^{\text{up}}_{\text{train}}$ &  $E_{\text{train}}$ &     $E^{\text{up}}_{\text{train}}$ &          $E_{\text{train}}$ &     $E^{\text{up}}_{\text{train}}$ &                        $E_{\text{train}}$ &    $E^{\text{up}}_{\text{train}}$ \\

\midrule
200                &  0.037 &    0.062 &    0.039 &      0.067 &       0.025 &         0.049 &  0.049 &    0.080 \\
400                &  0.024 &    0.044 &    0.027 &      0.047 &       0.017 &         0.030 &  0.031 &    0.059 \\
600                &  0.021 &    0.038 &    0.022 &      0.040 &       0.014 &         0.024 &  0.045 &    0.072 \\
800                &  0.018 &    0.032 &    0.019 &      0.035 &       0.013 &         0.020 &  0.031 &    0.056 \\
1000               &  0.017 &    0.029 &    0.018 &      0.032 &       0.013 &         0.018 &  0.030 &    0.054 \\
1200               &  0.016 &    0.028 &    0.017 &      0.029 &       0.012 &         0.015 &  0.029 &    0.052 \\
1400               &  0.015 &    0.026 &    0.016 &      0.028 &       0.012 &         0.014 &  0.028 &    0.051 \\
1600               &  0.015 &    0.025 &    0.015 &      0.027 &       0.012 &         0.013 &  0.029 &    0.052 \\
1800               &  0.014 &    0.024 &    0.015 &      0.026 &       0.011 &         0.012 &  0.030 &    0.053 \\
\bottomrule
\end{tabular}
\end{table}

\begin{table}[!th]
  \caption{Results for distinct $N_{\text{train}}$  averaged with respect to equations for 2D using loss \eqref{eq:upper_bound_loss}, $N_{\text{test}} = 200$.}
  \label{table:unsupervised_learning_for_different_test_size_for_all_models_averaged_wrt_equations}
  \centering
  \begin{tabular}{cllllllll}
\toprule
$N_{\text{train}}$  &                \multicolumn{2}{c}{fSNO} & \multicolumn{2}{c}{ChebNO} & \multicolumn{2}{c}{DilResNet} & \multicolumn{2}{c}{UNet} \\
{} &  $E_{\text{test}}$ &     $E^{\text{up}}_{\text{test}}$ &  $E_{\text{test}}$ &     $E^{\text{up}}_{\text{test}}$ &          $E_{\text{test}}$ &     $E^{\text{up}}_{\text{test}}$ &                        $E_{\text{test}}$ &    $E^{\text{up}}_{\text{test}}$ \\

\midrule
200  &  0.041  &  0.078  &  0.042 &  0.081 &  0.029 &  0.063 &         0.055 &  0.115 \\
400  &  0.027 &  0.057 & 0.029 &  0.059 &  0.019 &   0.039 &   0.038 &    0.080 \\
600   &  0.023 &    0.048 &  0.024 &  0.050 & 0.015 &  0.030 &   0.048 &    0.085 \\
800   &   0.020 &   0.042 &  0.021 &  0.044 &  0.014 &  0.025 &      0.034 &    0.068 \\
1000    &    0.018 &     0.038 & 0.020 &  0.040 &  0.013 &            0.022 &    0.033 &      0.066 \\
1200     &   0.018 &      0.036 &   0.018 &                               0.036 &     0.012 &    0.019 &     0.032 &    0.063 \\
1400     &      0.017  &   0.033 &   0.018 &   0.034 &   0.012 &  0.017 &     0.031 &   0.059 \\
1600    &     0.016 &    0.030 &  0.016 &  0.032 &  0.012 &            0.016 &    0.031 &      0.060 \\
1800  &   0.015 &   0.029 &    0.016 &  0.031 &   0.012 &  0.015 &     0.033 &     0.062 \\
\bottomrule
\end{tabular}
\end{table}

\begin{table}[!th]
  \caption{fSNO train and test results for all 2D datasets, $N_{\text{train}} = 1800$, $N_{\text{test}} = 200$.}
  \label{table:unsupervised_learning_fSNO_for_different_2D_datasets}
  \centering
   \begin{tabular}{lrrrrrr}
\toprule
equation  &  $E_{\text{train}}$ &  $E_{\text{test}}$ &  $E^{\text{up}}_{\text{train}}$ &  $E^{\text{up}}_{\text{test}}$ &  $R^{\text{up}}_{\text{train}}$ &  $R^{\text{up}}_{\text{test}}$ \\
\midrule
\eqref{eq:elliptic_1} &               0.012 &               0.013 &                           0.024 &                           0.031 &                            0.844 &                           0.821 \\
\eqref{eq:elliptic_3} &               0.017 &               0.018 &                           0.023 &                           0.023 &                            0.795 &                           0.819 \\
\eqref{eq:elliptic_4} &               0.014 &               0.017 &                           0.022 &                           0.028 &                            0.938 &                           0.929 \\
\eqref{eq:elliptic_6} &               0.013 &               0.014 &                           0.027 &                           0.034 &                            0.844 &                           0.834 \\
\bottomrule
\end{tabular}
\end{table}

\begin{table}[!th]
  \caption{ChebNO train and test results for all 2D datasets, $N_{\text{train}} = 1800$, $N_{\text{test}} = 200$.}
  \label{table:unsupervised_learning_ChebNO_for_different_2D_datasets}
  \centering
   \begin{tabular}{lrrrrrr}
\toprule
equation  &  $E_{\text{train}}$ &  $E_{\text{test}}$ &  $E^{\text{up}}_{\text{train}}$ &  $E^{\text{up}}_{\text{test}}$ &  $R^{\text{up}}_{\text{train}}$ &  $R^{\text{up}}_{\text{test}}$ \\
\midrule
\eqref{eq:elliptic_1} &               0.013 &               0.015 &                           0.028 &                           0.034 &                            0.858 &                           0.847 \\
\eqref{eq:elliptic_3} &               0.017 &               0.018 &                           0.025 &                           0.025 &                            0.767 &                           0.796 \\
\eqref{eq:elliptic_4} &               0.015 &               0.018 &                           0.024 &                           0.030 &                            0.941 &                           0.933 \\
\eqref{eq:elliptic_6} &               0.013 &               0.014 &                           0.028 &                           0.034 &                            0.859 &                           0.831 \\
\bottomrule
\end{tabular}
\end{table}

\begin{table}[!th]
  \caption{UNet train and test results for all 2D datasets, $N_{\text{train}} = 1800$, $N_{\text{test}} = 200$.}
  \label{table:unsupervised_learning_UNet_for_different_2D_datasets}
  \centering
   \begin{tabular}{lrrrrrr}
\toprule
equation  &  $E_{\text{train}}$ &  $E_{\text{test}}$ &  $E^{\text{up}}_{\text{train}}$ &  $E^{\text{up}}_{\text{test}}$ &  $R^{\text{up}}_{\text{train}}$ &  $R^{\text{up}}_{\text{test}}$ \\
\midrule
\eqref{eq:elliptic_1} &               0.012 &               0.014 &                           0.023 &                           0.035 &                            0.872 &                           0.823 \\
\eqref{eq:elliptic_3} &               0.025 &               0.024 &                           0.037 &                           0.036 &                            0.869 &                           0.883 \\
\eqref{eq:elliptic_4} &               0.070 &               0.078 &                           0.128 &                           0.142 &                            0.824 &                           0.819 \\
\eqref{eq:elliptic_6} &               0.013 &               0.014 &                           0.026 &                           0.035 &                            0.866 &                           0.820 \\
\bottomrule
\end{tabular}
\end{table}

\begin{table}[!th]
  \caption{DilResNet train and test results for all 2D datasets, $N_{\text{train}} = 1800$, $N_{\text{test}} = 200$.}
  \label{table:unsupervised_learning_DilResNet_for_different_2D_datasets}
  \centering
   \begin{tabular}{lrrrrrr}
\toprule
equation  &  $E_{\text{train}}$ &  $E_{\text{test}}$ &  $E^{\text{up}}_{\text{train}}$ &  $E^{\text{up}}_{\text{test}}$ &  $R^{\text{up}}_{\text{train}}$ &  $R^{\text{up}}_{\text{test}}$ \\
\midrule
\eqref{eq:elliptic_1}  &               0.008 &               0.009 &                           0.015 &                           0.018 &                            0.823 &                           0.771 \\
\eqref{eq:elliptic_3} &               0.010 &               0.012 &                           0.016 &                           0.017 &                            0.807 &                           0.826 \\
\eqref{eq:elliptic_4} &               0.009 &               0.011 &                           0.014 &                           0.019 &                            0.890 &                           0.881 \\
\eqref{eq:elliptic_6} &               0.009 &               0.010 &                           0.015 &                           0.018 &                            0.818 &                           0.803 \\
\bottomrule
\end{tabular}
\end{table}

\begin{table}[!th]
  \caption{Results for distinct equations averaged with respect to architectures, $N_{\text{train}}=1200$.}
  \label{table:eq_2D}
  \centering
  \begin{tabular}{ccccccccccc}
    \toprule
     equation & $E_{\text{train}}$ & $E_{\text{test}}$ & $E^{\text{ub}}_{\text{train}}$ & $E^{\text{ub}}_{\text{test}}$ & $R^{\text{ub}}_{\text{train}}$ & $R^{\text{ub}}_{\text{test}}$ & $\widetilde{E}^{\text{ub}}_{\text{train}}$ & $\widetilde{E}^{\text{ub}}_{\text{test}}$ & $\widetilde{R}^{\text{ub}}_{\text{train}}$ & $\widetilde{R}^{\text{ub}}_{\text{test}}$ \\
    \cmidrule(r){1-11}
    \eqref{eq:elliptic_1} & $0.009$ & $0.012$ & $0.017$ & $0.023$ & $0.861$ & $0.859$ & $0.117$ & $0.125$ & $0.840$ & $0.831$ \\
    \eqref{eq:elliptic_3} & $0.015$ & $0.017$ & $0.031$ & $0.034$ & $0.980$ & $0.986$ & $0.103$ & $0.107$ & $0.832$ & $0.841$ \\
    \eqref{eq:elliptic_4} & $0.015$ & $0.020$ & $0.027$ & $0.034$ & $0.977$ & $0.984$ & $0.113$ & $0.138$ & $0.752$ & $0.794$ \\
    \eqref{eq:elliptic_6} & $0.010$ & $0.013$ & $0.017$ & $0.021$ & $0.861$ & $0.877$ & $0.111$ & $0.116$ & $0.789$ & $0.786$ \\
    \bottomrule
\end{tabular}
\end{table}

\begin{table}[!th]
  \caption{Results for distinct architectures averaged with respect to equations, $N_{\text{train}}=1200$.}
  \label{table:model_2D}
  \centering
  \begin{tabular}{ccccccccccc}
    \toprule
     Architecture & $E_{\text{train}}$ & $E_{\text{test}}$ & $E^{\text{ub}}_{\text{train}}$ & $E^{\text{ub}}_{\text{test}}$ & $R^{\text{ub}}_{\text{train}}$ & $R^{\text{ub}}_{\text{test}}$ & $\widetilde{E}^{\text{ub}}_{\text{train}}$ & $\widetilde{E}^{\text{ub}}_{\text{test}}$ & $\widetilde{R}^{\text{ub}}_{\text{train}}$ & $\widetilde{R}^{\text{ub}}_{\text{test}}$ \\
    \cmidrule(r){1-11}
    ChebNO & $0.012$ & $0.014$ & $0.021$ & $0.024$ & $0.916$ & $0.922$ & $0.092$ & $0.097$ & $0.911$ & $0.906$ \\
    DilResNet & $0.013$ & $0.016$ & $0.023$ & $0.035$ & $0.939$ & $0.935$ & $0.072$ & $0.077$ & $0.779$ & $0.838$ \\
    FNO & $0.009$ & $0.016$ & $0.017$ & $0.026$ & $0.911$ & $0.935$ & $0.133$ & $0.168$ & $0.818$ & $0.838$ \\
    MLP & $0.020$ & $0.020$ & $0.036$ & $0.034$ & $0.922$ & $0.916$ & $0.166$ & $0.167$ & $0.709$ & $0.693$ \\
    UNet & $0.009$ & $0.013$ & $0.018$ & $0.023$ & $0.923$ & $0.929$ & $0.090$ & $0.096$ & $0.783$ & $0.777$ \\
    fSNO & $0.011$ & $0.013$ & $0.022$ & $0.025$ & $0.907$ & $0.921$ & $0.114$ & $0.124$ & $0.820$ & $0.827$ \\
    \bottomrule
\end{tabular}
\end{table}

\begin{table}[!th]
  \caption{Results for distinct $N_{\text{train}}$ averaged with respect to the network type and equation.}
  \label{table:N_train_2D}
  \centering
  \begin{tabular}{ccccccccccc}
    \toprule
     $N_{\text{train}}$ & $E_{\text{train}}$ & $E_{\text{test}}$ & $E^{\text{ub}}_{\text{train}}$ & $E^{\text{ub}}_{\text{test}}$ & $R^{\text{ub}}_{\text{train}}$ & $R^{\text{ub}}_{\text{test}}$ & $\widetilde{E}^{\text{ub}}_{\text{train}}$ & $\widetilde{E}^{\text{ub}}_{\text{test}}$ & $\widetilde{R}^{\text{ub}}_{\text{train}}$ & $\widetilde{R}^{\text{ub}}_{\text{test}}$ \\
    \cmidrule(r){1-11}
200 & $0.039$ & $0.052$ & $0.113$ & $0.190$ & $0.916$ & $0.925$ & $0.263$ & $0.322$ & $0.768$ & $0.819$ \\
400 & $0.023$ & $0.030$ & $0.049$ & $0.073$ & $0.917$ & $0.930$ & $0.177$ & $0.206$ & $0.795$ & $0.813$ \\
600 & $0.018$ & $0.023$ & $0.033$ & $0.047$ & $0.915$ & $0.925$ & $0.145$ & $0.165$ & $0.794$ & $0.809$ \\
800 & $0.015$ & $0.019$ & $0.027$ & $0.035$ & $0.915$ & $0.925$ & $0.135$ & $0.149$ & $0.804$ & $0.818$ \\
1000 & $0.013$ & $0.017$ & $0.024$ & $0.031$ & $0.914$ & $0.924$ & $0.116$ & $0.129$ & $0.795$ & $0.809$ \\
1200 & $0.012$ & $0.015$ & $0.023$ & $0.028$ & $0.920$ & $0.926$ & $0.111$ & $0.122$ & $0.803$ & $0.813$ \\
\bottomrule
\end{tabular}
\end{table}

\begin{table}[!th]
  \caption{Equation \eqref{eq:elliptic_1}, different resolutions, results are averaged with respect to architecture type.}
  \label{table:J_2D}
  \centering
  \begin{tabular}{ccccccccc}
    \toprule
     $J$ & $E_{\text{train}}$ & $E_{\text{test}}$ & $E^{\text{ub}}_{\text{train}}$ & $E^{\text{ub}}_{\text{test}}$ & $\widetilde{E}^{\text{ub}}_{\text{train}}$ & $\widetilde{E}^{\text{ub}}_{\text{test}}$ & $\widetilde{E}^{\text{ub}}_{\text{train}} \big/ E^{\text{ub}}_{\text{train}}$ & $\widetilde{E}^{\text{ub}}_{\text{test}} \big/ E^{\text{ub}}_{\text{test}}$ \\
    \cmidrule(r){1-9}
    $5$ & $0.009$ & $0.012$ & $0.017$ & $0.023$ & $0.117$ & $0.125$ & $7.019$ & $5.448$ \\
    $6$ & $0.012$ & $0.015$ & $0.024$ & $0.027$ & $0.127$ & $0.135$ & $5.364$ & $5.020$ \\
    $7$ & $0.020$ & $0.022$ & $0.060$ & $0.058$ & $0.183$ & $0.183$ & $3.043$ & $3.174$ \\
    \bottomrule
\end{tabular}

\end{table}

\begin{table}[!th]
  \caption{Ratios $\frac{\text{upper bound obtained by a posteriori neural solver}}{\text{upper bound obtained by optimization}}$, $N_{\text{train}}=1200$.}
  \label{table:quality_2D}
  \centering
  \begin{tabular}{ccccccccc}
\toprule
 & \multicolumn{4}{c}{train set} & \multicolumn{4}{c}{test set} \\
architecture\,\textbackslash\,equation & \eqref{eq:elliptic_1} & \eqref{eq:elliptic_3} & \eqref{eq:elliptic_4} & \eqref{eq:elliptic_6} & \eqref{eq:elliptic_1} & \eqref{eq:elliptic_3} & \eqref{eq:elliptic_4} & \eqref{eq:elliptic_6} \\
\cmidrule(r){1-9}
ChebNO & $4.5$ & $3.8$ & $3.0$ & $5.2$ & $5.4$ & $4.0$ & $3.3$ & $5.9$ \\
DilResNet & $2.7$ & $1.5$ & $2.2$ & $3.3$ & $4.0$ & $2.4$ & $2.9$ & $3.6$ \\
FNO & $6.5$ & $5.9$ & $6.6$ & $7.3$ & $10.1$ & $5.6$ & $6.7$ & $10.0$ \\
MLP & $8.6$ & $1.6$ & $5.7$ & $5.8$ & $8.9$ & $1.6$ & $5.0$ & $6.0$ \\
UNet & $4.6$ & $3.3$ & $3.6$ & $5.8$ & $6.9$ & $3.4$ & $3.4$ & $8.4$ \\
fSNO & $5.2$ & $4.8$ & $5.0$ & $5.5$ & $5.9$ & $4.9$ & $4.8$ & $6.1$ \\
\bottomrule
\end{tabular}
\end{table}

\begin{table}[!th]
  \caption{Results for fSNO for different  $\lambda$  for 2D datasets, $N_{\text{train}}=1200$, $N_{\text{test}}=800$}.
  \label{table:unsupervised_learning_fSNO_2D_for_different_lambda}
  \centering
  \begin{tabular}{lllllllllll}
\toprule
equation &              \multicolumn{2}{c}{$\lambda = 0.01$} &                                \multicolumn{2}{c}{$\lambda = 0.1$} &                                \multicolumn{2}{c}{$\lambda = 1$} &                                   \multicolumn{2}{c}{$\lambda = 10$} &                             \multicolumn{2}{c}{$\lambda = 100$}                \\
                   &  $E_{\text{train}}$ &     $E_{\text{test}}$ &    $E_{\text{train}}$ &     $E_{\text{test}}$ &  $E_{\text{train}}$ &     $E_{\text{test}}$ &    $E_{\text{train}}$ &  $E_{\text{test}}$ &  $E_{\text{train}}$ &     $E_{\text{test}}$ \\
\midrule
\eqref{eq:elliptic_1} &  0.026 &               0.027 &                0.016 &  0.018 &                0.014 &                 0.016 &               0.017 &               0.019 &  0.104 &  0.102 \\
\eqref{eq:elliptic_3} &                 0.028 &               0.028 &                 0.02 &                  0.02 &                0.019 &                 0.019 &               0.021 &               0.021 &                 0.03 &                0.029 \\
\eqref{eq:elliptic_4} &                 0.021 &               0.023 &                0.015 &                 0.017 &                0.016 &  0.018 &               0.017 &               0.019 &  0.071 &                 0.070 \\
\eqref{eq:elliptic_6} &                 0.032 &               0.033 &                0.016 &                 0.019 &                0.014 &                 0.017 &               0.017 &               0.019 &                0.105 &                0.105 \\
\bottomrule
\end{tabular}
\end{table}

\begin{table}[!th]
  \caption{Results for ChebNO for different  $\lambda$  for 2D datasets, $N_{\text{train}}=1200$, $N_{\text{test}}=800$}.
  \label{table:unsupervised_learning_ChebNO_2D_for_different_lambda}
  \centering
  \begin{tabular}{lllllllllll}
\toprule
equation &              \multicolumn{2}{c}{$\lambda = 0.01$} &                                \multicolumn{2}{c}{$\lambda = 0.1$} &                                \multicolumn{2}{c}{$\lambda = 1$} &                                   \multicolumn{2}{c}{$\lambda = 10$} &                             \multicolumn{2}{c}{$\lambda = 100$}                \\
                   &  $E_{\text{train}}$ &     $E_{\text{test}}$ &    $E_{\text{train}}$ &     $E_{\text{test}}$ &  $E_{\text{train}}$ &     $E_{\text{test}}$ &    $E_{\text{train}}$ &  $E_{\text{test}}$ &  $E_{\text{train}}$ &     $E_{\text{test}}$ \\
\midrule
\eqref{eq:elliptic_1} &  0.026 &                 0.028 &                0.016 &  0.018 &                0.015 &                 0.017 &               0.016 &  0.018 &                0.064 &                 0.062 \\
\eqref{eq:elliptic_3} &                 0.032 &                 0.031 &                0.021 &                 0.021 &                 0.02 &                  0.02 &               0.021 &                 0.021 &                0.028 &                 0.027 \\
\eqref{eq:elliptic_4} &                 0.024 &  0.026 &                0.016 &  0.018 &                0.016 &  0.018 &                0.02 &                 0.021 &                0.055 &                 0.054 \\
\eqref{eq:elliptic_6} &                  0.03 &                 0.031 &                0.017 &                 0.019 &                0.016 &                 0.017 &               0.019 &                 0.021 &                0.025 &  0.026 \\
\bottomrule
\end{tabular}
\end{table}

\begin{table}[!th]
  \caption{Results for UNet for different  $\lambda$  for 2D datasets, $N_{\text{train}}=1200$, $N_{\text{test}}=800$}.
  \label{table:unsupervised_learning_UNet_2D_for_different_lambda}
  \centering
  \begin{tabular}{lllllllllll}
\toprule
equation &              \multicolumn{2}{c}{$\lambda = 0.01$} &                                \multicolumn{2}{c}{$\lambda = 0.1$} &                                \multicolumn{2}{c}{$\lambda = 1$} &                                   \multicolumn{2}{c}{$\lambda = 10$} &                             \multicolumn{2}{c}{$\lambda = 100$}                \\
                   &  $E_{\text{train}}$ &     $E_{\text{test}}$ &    $E_{\text{train}}$ &     $E_{\text{test}}$ &  $E_{\text{train}}$ &     $E_{\text{test}}$ &    $E_{\text{train}}$ &  $E_{\text{test}}$ &  $E_{\text{train}}$ &     $E_{\text{test}}$ \\
\midrule
\eqref{eq:elliptic_1} &                 0.02 &               0.021 &  0.013 &                 0.015 &                0.015 &               0.015 &                0.021 &               0.022 &  0.104 &  0.102 \\
\eqref{eq:elliptic_3} &                 0.04 &                0.04 &                 0.028 &                 0.027 &                0.019 &               0.019 &                0.091 &                0.09 &                 0.11 &                0.109 \\
\eqref{eq:elliptic_4} &                0.024 &               0.024 &   0.071 &                  0.07 &  0.071 &                0.07 &  0.071 &                0.07 &  0.071 &                 0.07 \\
\eqref{eq:elliptic_6} &                0.016 &               0.017 &                 0.016 &  0.018 &                0.012 &               0.015 &                0.032 &               0.035 &                0.105 &                0.105 \\
\bottomrule
\end{tabular}
\end{table}

\begin{table}[!th]
  \caption{Results for DilResNet for different  $\lambda$  for 2D datasets, $N_{\text{train}}=1200$, $N_{\text{test}}=800$}.
  \label{table:unsupervised_learning_DilResNet_2D_for_different_lambda}
  \centering
  \begin{tabular}{lllllllllll}
\toprule
equation &              \multicolumn{2}{c}{$\lambda = 0.01$} &                                \multicolumn{2}{c}{$\lambda = 0.1$} &                                \multicolumn{2}{c}{$\lambda = 1$} &                                   \multicolumn{2}{c}{$\lambda = 10$} &                             \multicolumn{2}{c}{$\lambda = 100$}                \\
                   &  $E_{\text{train}}$ &     $E_{\text{test}}$ &    $E_{\text{train}}$ &     $E_{\text{test}}$ &  $E_{\text{train}}$ &     $E_{\text{test}}$ &    $E_{\text{train}}$ &  $E_{\text{test}}$ &  $E_{\text{train}}$ &     $E_{\text{test}}$ \\
\midrule
\eqref{eq:elliptic_1} &  0.013 &  0.013 &                  0.01 &                  0.01 &  0.009 &                0.01 &  0.009 &                  0.01 &                 0.017 &                 0.017 \\
\eqref{eq:elliptic_3} &                 0.021 &                  0.02 &  0.018 &  0.018 &  0.018 &               0.018 &  0.018 &  0.018 &                 0.027 &  0.026 \\
\eqref{eq:elliptic_4} &  0.013 &                 0.014 &                 0.012 &                 0.012 &                 0.011 &               0.012 &                 0.012 &                 0.012 &                 0.048 &                 0.047 \\
\eqref{eq:elliptic_6} &  0.013 &                 0.014 &  0.009 &                  0.01 &                  0.01 &                0.01 &                  0.01 &                 0.011 &  0.026 &  0.026 \\
\bottomrule
\end{tabular}
\end{table}
\end{document}


\maketitle

\section{Code overview}
We use JAX \cite{jax2018github}, Equinox \cite{kidger2021equinox} and Optax \cite{deepmind2020jax}.

The GitHub repository with the code will be made available after the review stage. During the review we provide two GoogleColab notebooks that are sufficient to reproduce a portion of results from the article.

Supervised: \href{https://colab.research.google.com/drive/18V2CYcepCAWC4siNkyxEcUyjskm7EcWs}{https://colab.research.google.com/drive/18V2CYcepCAWC4siNkyxEcUyjskm7EcWs}

Unsupervised: \href{https://colab.research.google.com/drive/1OI-jFycmpE2kPt8PJ4H9dMxoU7yL9NeE}{https://colab.research.google.com/drive/1OI-jFycmpE2kPt8PJ4H9dMxoU7yL9NeE}

\section{Details on architectures used}
\subsection{Supervised training}
Here we provide more details on the architectures used.
\begin{enumerate}
    \item Construction of FNO closely follows the one given in \cite{li2020fourier}. We use $24$ featrures in processor, $4$ layers and the number of modes is $\left\lceil N_{\text{spatial}} / 4\right\rceil$.
    \item Construction of fSNO is similar to the one of FNO (i.e., encoder-processor-decoder and integral kernel in place of the linear layer), but following \cite{fanaskov2022spectral} we replace Fourier basis with orthogonal polynomials in the integral kernel. For this particular architecture, the construction of the integral kernel is as follows. We use the Gauss-Chebyshev grid, and compute projection on the basis of polynomials using Gauss quadratures \cite{golub1969calculation}. After the projection on the space of polynomials, we apply three convolutions with kernel size $3$. Finally, to return to the physical space, we compute sum $\sum_{n} c_n p_n(x)$ on the Gauss-Chebyshev grid, where $c_n$ are coefficients obtained after convolutions. In this case, the number of features in the processor is $34$, the number of layers is $4$ and the number of modes is $\left\lceil N_{\text{spatial}} / 4\right\rceil$
    \item ChebNO is a spectral neural operator \cite{fanaskov2022spectral} defined on Chebyshev grid. Construction of integral kernel is similar to fSNO but DCT is used to find projection on the polynomial space in place of Gauss quadratures. In the current architecture, the processor has $32$ features, the number of layers is $4$, and $16$ modes are used for all grids.
    \item DilResNet that we use closely follows architecture described in \cite{stachenfeld2021learned}. Namely, we use a processor with $24$ features that has $4$ layers. Each layer consists of convolutions with strides $[1, 2, 4, 8, 4, 2, 1]$, kernel size $3$, and skip connection.
    \item For UNet \cite{ronneberger2015u} we start with $10$ features and double the number of features with each downsampling that decreases the number of grid points by the factor of $2$ in each dimension. On each grid we use $2$ convolutions (kernel size $3$) and max pooling, transposed convolution are used for upsampling, and $3$ convolutions (kernel size $3$) appears on each grid after upsampling. In total, we have $4$ grids.
    \item MLP that we use consists of linear layers that process each dimension (including the feature dimension) separately. That way, the linear operator is defined by three matrices in $D=2$. MLP uses $64$ in the processor and has $4$ layers.
\end{enumerate}
For all networks, we use $\text{ReLU}$ nonlinearity.
Typical number of parameters for each network is given in the table below.

\begin{center}
\begin{tabular}{lcccccc}
\toprule
& FNO & fSNO & ChebNO & DilResNet & UNet & MLP \\
\# parameters& $668\times 10^{3}$ & $130\times 10^{3}$ & $115\times 10^{3}$ & $147\times 10^{3}$ & $248\times 10^3$ & $24\times 10^3$\\
\bottomrule
\end{tabular}
\end{center}

\subsection{Unsupervised training}

Here we provide more details on the architectures used for 1D and 2D datasets. The architectures both for 1D and 2D cases are the 1D  and 2D versions of the architectures described previously.

For 1D datasets, the detailed description of the architectures is given below.
\begin{enumerate}
    \item For fSNO, after the projection on the space of polynomials, we apply three convolutions with kernel size $3$. In this case, the number of features in the processor is $34$, the number of layers is $4$ and the number of modes is $\left\lceil N_{\text{spatial}} / 4\right\rceil$.
    \item For ChebNO, the architecture has processor with $32$ features, $4$ layers and $16$ modes.
    \item For DilResNet, we use a processor with $64$ features that has $5$ layers. Each layer consists of convolutions with strides $[1, 2, 4, 8, 4, 2, 1]$, kernel size $3$, and skip connection.
    \item For UNet, we start with $32$ features. Kernel size in all convolutions is $3$. In total, we have $4$ grids.
\end{enumerate}

For 2D datasets, we provide the detailed description of the architectures below.
\begin{enumerate}
    \item For fSNO, kernel size is $3$, the number of features in the processor is $32$, the number of layers is $4$ and the number of modes is $\left\lceil N_{\text{spatial}} / 4\right\rceil$.
    \item For ChebNO, the architecture is similar to supervised case.
    \item For DilResNet, a processor has $32$ features and $5$ layers. Each layer consists of convolutions with strides $[1, 2, 4, 8, 4, 2, 1]$, kernel size $3$, and skip connection.
    \item For UNet, we start with $10$ features. On each grid we use $2$ convolutions with kernel size equals $3$ and max pooling. In total, we have $4$ grids.
\end{enumerate}

We use ReLU nonlinear activation for all architectures for 1D and 2D cases. The number of parameters for all neural operators are given below.

\begin{center}
\begin{tabular}{lcccc}
\toprule
& \multicolumn{4}{c}{\# parameters} \\
& fSNO & ChebNO & DilResNet & UNet  \\
\midrule
$D=1$ & $47\times 10^{3}$ & $12 \times 10^{3}$ &$433\times 10^{3}$ & $19\times 10^{5}$\\
$D=2$ & $115 \times 10^{3}$ & $116 \times 10^{3}$ &$326 \times 10^{3}$ & $248 \times 10^{3}$\\
\bottomrule
\end{tabular}
\end{center}

\section{Additional results}
\subsection{Supervised training}
For convenience we recall parameters of $D=2$ elliptic PDE
\begin{equation}
    \label{eq:general_elliptic_equation}
    \begin{split}
    -\sum_{ij}\frac{\partial}{\partial x_i}\left(a_{ij}(x) \frac{\partial u}{\partial x_j} \right) + b^{2}(x) u(x) = f(x),\,x\in\Gamma = [0, 1]^{D},\,\left.u\right|_{\partial\Gamma} = 0,\,a_{ij}(x) \geq c > 0,
    \end{split}
\end{equation}
used in experiments.

Random polynomials that we use read
\begin{equation}
    \label{eq:random_trig_2D}
    \mathcal{P}(N_1, N_2, \alpha) = \left\{f(x) = \mathcal{R}\left(\sum_{m=0}^{N_1}\sum_{n=0}^{N_2}\frac{c_{mn}\exp\left(2\pi i(mx_1 + nx_2)\right)}{(1+m+n)^\alpha}\right):\mathcal{R}(c), \mathcal{I}(c)\simeq \mathcal{N}(0, I)\right\}.
\end{equation}

With these polynomials we define $4$ test equations
\begin{equation}
    \label{eq:elliptic_1}
    \begin{split}
    &a(x) = \begin{pmatrix}
        \alpha(x) & 0 \\
        \gamma(x) & \beta(x)
    \end{pmatrix}
    \begin{pmatrix}
        \alpha(x) & \gamma(x) \\
        0 & \beta(x)
    \end{pmatrix},\\
    &\,\alpha(x),\beta(x) \simeq 0.1\mathcal{P}(5, 5, 2) + 1;\,\gamma(x),\,b(x),\,f(x) \simeq \mathcal{P}(5, 5, 2).
    \end{split}
\end{equation}

\begin{equation}
    \label{eq:elliptic_3}
    \begin{split}
    a(x) = \alpha(x)I,\,\alpha(x) = \begin{cases}
        10,\,p_1(x) \geq 0;\\
        1,\, p_1(x) < 0,
    \end{cases}\,b(x) = 0,\,f(x) = 1,\,p_1(x) \simeq \mathcal{P}(5, 5, 2).
    \end{split}
\end{equation}

\begin{equation}
    \label{eq:elliptic_4}
    \begin{split}
    a(x) = \alpha(x)I,\,\alpha(x) = \begin{cases}
        10,\,p_1(x) \geq 0;\\
        1,\, p_1(x) < 0,
    \end{cases}\,b(x),\,f(x),\,p_1(x) \simeq \mathcal{P}(5, 5, 2).
    \end{split}
\end{equation}

\begin{equation}
    \label{eq:elliptic_6}
    \begin{split}
    &a(x) = \begin{pmatrix}
        \alpha(x) & 0 \\
        \gamma(x) & \beta(x)
    \end{pmatrix}
    \begin{pmatrix}
        \alpha(x) & \gamma(x) \\
        0 & \beta(x)
    \end{pmatrix},\\
    &\,\alpha(x),\beta(x) \simeq 0.1\mathcal{P}(5, 5, 2) + 1;\,\gamma(x), f(x) \simeq \mathcal{P}(5, 5, 2);\,b(x) = 0.
    \end{split}
\end{equation}

In the Tables \ref{table:ChebNO elliptic ub}, \ref{table:DilResNet elliptic ub}, \ref{table:FNO elliptic ub}, \ref{table:fSNO elliptic ub}, \ref{table:MLP elliptic ub}, \ref{table:UNet elliptic ub} one can find results for train runs with different train set sizes $N_{\text{train}}$ for all architectures and equations used in the main text (number of points along each dimension is $2^5+1$).

\subsection{Unsupervised training}
For the unsupervised case we consider $D=1$ version of \eqref{eq:general_elliptic_equation}
and optimize upper bound to find approximate solution $v$ and upper bound certificate $y$:
\begin{equation}
    \label{eq:upper_bound_1D}
     \left\|v - u\right\|_{a}  \leq\inf_{v, y}  \left(\left(\int_{\Gamma} \frac{1}{a}\left( y - a\frac{d v}{dx} \right)^2 dx\right)^{1/2} + \overline{C}_{F} \left\| f + \frac{d y}{dx} - b v\right\|_2\right),\,\overline{C}_{F} = \frac{\sup a(x)}{\pi}.
\end{equation}
To draw random functions, we use the space analogous to \eqref{eq:random_trig_2D} and denote it $\mathcal{P}(N, \alpha)$, where $N$ controls the number of modes and $\alpha$ controls the decay of coefficients.

The parameters for the first equation reads
\begin{equation}
    \label{eq:elliptic_1_1D}
    a(x) \sim 0.1\mathcal{P}(5, 0)^2 + 1;\,b(x) \sim 0.2\mathcal{P}(5, 0),\,f(x)\sim\mathcal{P}(5, 0).
\end{equation}
The second equation is similar but with zero source terms, i.e.,
\begin{equation}
    \label{eq:elliptic_2_1D}
    a(x) \sim 0.1\mathcal{P}(5, 0)^2 + 1;\,b(x)=0,\,f(x)\sim\mathcal{P}(5, 0).
\end{equation}
The third equation has discontinuous diffusion coefficient $a$, shared randomly selected right-hand side $f_1(x) \sim \mathcal{P}(5, 0)$ for all samples and zero sources term $b(x) = 0$. The diffusion coefficient reads
\begin{equation}
    \label{eq:elliptic_3_1D}
    a(x) = \begin{cases}
        10,\, x \in [\alpha, \beta];\\
        1,\, x \notin [\alpha, \beta],
        \\
    \end{cases}\alpha=\min(A, B),\,\beta=\max(A,B),\,A,B\sim U([0, 1]),
\end{equation}
where $U([0, 1])$ is a uniform distribution on the interval $[0, 1]$.

For training neural network we used the following loss function:
\[
\mathcal{L}\left(y, v, a, b, f, \overline{C}_{F}\right) = \sqrt{\int_{\Gamma} \frac{1}{a}\left( y - a\frac{d v}{dx} \right)^2 dx} + \overline{C}_{F} \left\| f + \frac{d y}{dx} - b v\right\|_2 + \lambda \sqrt{v^2_{\partial \Gamma}},
\]
\[
\mathcal{L}\left(y, v, a, b, f, \overline{C}_{F}\right) \rightarrow \min_{v,~y},
\]
where $y, v$ are the output of the neural network,  $a, b, f$ are defined for each equations \eqref{eq:elliptic_1_1D}, \eqref{eq:elliptic_2_1D}, \eqref{eq:elliptic_3_1D} and $\lambda$ is the hyperparameter of the loss function.  The part of the loss function with boundaries is required as neural network cannot learn the values of the exact solution for $x \in \Gamma$ as upper bound does not contain such information.




For unsupervised training in $D=1$ datasets, we train for $500$ epochs with Adam optimizer that has the learning rate $5 \cdot 10^{-3}$ multiplied by $0.5$ for each $50$ epoch and weight decay $10^{-2}$. Training parameters for 2D case are described in the main text.

In Tables \ref{table:fsno_lambda}, \ref{table:chebno_lambda}, \ref{table:Dilresnet_lambda}, \ref{table:u_lambda} one can find results for train runs with different $\lambda$ for all architectures and equations for 1D datasets (the number of points along each dimension is $101$). In  the Tables \ref{table:fsno_datasets}, \ref{table:chebno_datasets}, \ref{table:Dilresnet_datasets} and \ref{table:u_datasets}, the detailed results for all models and equations with optimal $\lambda$ ($\lambda = 1$ ) are presented. We provide an aggregated summary of the results for the $1D$ datasets in the Tables~\ref{table:unsup_1D_eq},\,\ref{table:unsup_1D_m},\,\ref{table:unsup_1D_trsize},\,\ref{table:unsup_1D_unet}. For each network and PDE we follow the same experiment protocol and train the network with $N_{\text{train}} = 200, 400, 600, 800, 1000, 1200$. In all cases, we use $N=101$ grid points.

In the Tables \ref{table:fsno_lambda_2D}, \ref{table:chebno_lambda_2D}, \ref{table:dilresnet_lambda_2D}, \ref{table:unet_lambda_2D} one can find results for train runs with different $\lambda$ for all architectures and equations for 2D datasets. In the Tables \ref{table:fsno_datasets_2d}, \ref{table:chebno_datasets_2d}, \ref{table:dilresnet_datasets_2d}, \ref{table:unet_datasets_2d}, the detailed results for all models and equations are presented. An aggregated summary of the results is given in the main text.

Overall, the trends for 1D and 2D  are mostly the same: the increase in the number of train points decreases the test error, the correlations between errors and upper bounds are good even when the bounds are not tight.





\begin{table}[!ht]
\centering
\caption{ChebNO}
\label{table:ChebNO elliptic ub}
\resizebox{\textwidth}{!}{
\begin{tabular}{llrrrrrrrrrr}
\toprule
 equation & $N_{\text{train}}$  & $E_{\text{train}}$ & $E_{\text{test}}$ & $E^{\text{ub}}_{\text{train}}$ & $E^{\text{ub}}_{\text{test}}$ & $R^{\text{ub}}_{\text{train}}$ & $R^{\text{ub}}_{\text{test}}$ & $\widetilde{E}^{\text{ub}}_{\text{train}}$ & $\widetilde{E}^{\text{ub}}_{\text{test}}$ & $\widetilde{R}^{\text{ub}}_{\text{train}}$ & $\widetilde{R}^{\text{ub}}_{\text{test}}$ \\
\midrule
\multirow[t]{6}{*}{\eqref{eq:elliptic_1}, $J=5$} & 200 & 0.023 & 0.029 & 0.041 & 0.054 & 0.861 & 0.831 & 0.121 & 0.143 & 0.840 & 0.833 \\
 & 400 & 0.016 & 0.019 & 0.028 & 0.038 & 0.839 & 0.852 & 0.106 & 0.117 & 0.871 & 0.887 \\
 & 600 & 0.013 & 0.016 & 0.023 & 0.031 & 0.831 & 0.849 & 0.098 & 0.105 & 0.904 & 0.913 \\
 & 800 & 0.011 & 0.013 & 0.019 & 0.026 & 0.839 & 0.847 & 0.095 & 0.099 & 0.906 & 0.904 \\
 & 1000 & 0.010 & 0.012 & 0.017 & 0.022 & 0.831 & 0.838 & 0.089 & 0.093 & 0.905 & 0.898 \\
 & 1200 & 0.009 & 0.011 & 0.017 & 0.020 & 0.851 & 0.844 & 0.091 & 0.091 & 0.908 & 0.892 \\
\cline{1-12}
\multirow[t]{6}{*}{\eqref{eq:elliptic_1}, $J=6$} & 200 & 0.033 & 0.040 & 0.065 & 0.073 & 0.880 & 0.881 & 0.126 & 0.144 & 0.814 & 0.823 \\
\bottomrule
\end{tabular}
}

\end{table}

\begin{table}[!ht]
\centering
\caption{DilResNet}
\label{table:DilResNet elliptic ub}
\resizebox{\textwidth}{!}{
\begin{tabular}{llrrrrrrrrrr}
\toprule
 equation & $N_{\text{train}}$  & $E_{\text{train}}$ & $E_{\text{test}}$ & $E^{\text{ub}}_{\text{train}}$ & $E^{\text{ub}}_{\text{test}}$ & $R^{\text{ub}}_{\text{train}}$ & $R^{\text{ub}}_{\text{test}}$ & $\widetilde{E}^{\text{ub}}_{\text{train}}$ & $\widetilde{E}^{\text{ub}}_{\text{test}}$ & $\widetilde{R}^{\text{ub}}_{\text{train}}$ & $\widetilde{R}^{\text{ub}}_{\text{test}}$ \\
\midrule
\multirow[t]{6}{*}{\eqref{eq:elliptic_1}, $J=5$} & 200 & 0.050 & 0.073 & 0.126 & 0.504 & 0.882 & 0.954 & 0.236 & 0.365 & 0.776 & 0.943 \\
 & 400 & 0.029 & 0.037 & 0.056 & 0.114 & 0.933 & 0.955 & 0.131 & 0.163 & 0.828 & 0.933 \\
 & 600 & 0.019 & 0.024 & 0.035 & 0.082 & 0.922 & 0.942 & 0.090 & 0.112 & 0.868 & 0.939 \\
 & 800 & 0.014 & 0.017 & 0.026 & 0.043 & 0.923 & 0.942 & 0.076 & 0.090 & 0.884 & 0.949 \\
 & 1000 & 0.010 & 0.013 & 0.019 & 0.036 & 0.912 & 0.928 & 0.063 & 0.072 & 0.880 & 0.932 \\
 & 1200 & 0.008 & 0.011 & 0.016 & 0.025 & 0.915 & 0.920 & 0.062 & 0.069 & 0.884 & 0.927 \\
\cline{1-12}
\multirow[t]{6}{*}{\eqref{eq:elliptic_1}, $J=6$} & 200 & 0.086 & 0.108 & 0.488 & 0.716 & 0.872 & 0.883 & 0.621 & 0.624 & 0.842 & 0.778 \\
\bottomrule
\end{tabular}
}
\end{table}

\begin{table}[!ht]
\centering
\caption{FNO}
\label{table:FNO elliptic ub}
\resizebox{\textwidth}{!}{
\begin{tabular}{llrrrrrrrrrr}
\toprule
 equation & $N_{\text{train}}$  & $E_{\text{train}}$ & $E_{\text{test}}$ & $E^{\text{ub}}_{\text{train}}$ & $E^{\text{ub}}_{\text{test}}$ & $R^{\text{ub}}_{\text{train}}$ & $R^{\text{ub}}_{\text{test}}$ & $\widetilde{E}^{\text{ub}}_{\text{train}}$ & $\widetilde{E}^{\text{ub}}_{\text{test}}$ & $\widetilde{R}^{\text{ub}}_{\text{train}}$ & $\widetilde{R}^{\text{ub}}_{\text{test}}$ \\
\midrule
\multirow[t]{6}{*}{\eqref{eq:elliptic_1}, $J=5$} & 200 & 0.011 & 0.036 & 0.019 & 0.065 & 0.852 & 0.882 & 0.312 & 0.385 & 0.629 & 0.881 \\
 & 400 & 0.010 & 0.025 & 0.017 & 0.050 & 0.829 & 0.884 & 0.195 & 0.246 & 0.804 & 0.807 \\
 & 600 & 0.009 & 0.021 & 0.015 & 0.036 & 0.848 & 0.871 & 0.171 & 0.221 & 0.844 & 0.784 \\
 & 800 & 0.009 & 0.018 & 0.015 & 0.032 & 0.837 & 0.884 & 0.145 & 0.200 & 0.822 & 0.833 \\
 & 1000 & 0.008 & 0.017 & 0.014 & 0.030 & 0.840 & 0.879 & 0.144 & 0.177 & 0.815 & 0.829 \\
 & 1200 & 0.008 & 0.015 & 0.014 & 0.026 & 0.842 & 0.877 & 0.137 & 0.170 & 0.836 & 0.825 \\
\cline{1-12}
\multirow[t]{6}{*}{\eqref{eq:elliptic_1}, $J=6$} & 200 & 0.014 & 0.044 & 0.028 & 0.077 & 0.892 & 0.847 & 0.265 & 0.429 & 0.654 & 0.794 \\
\bottomrule
\end{tabular}
}
\end{table}

\begin{table}[!ht]
\centering
\caption{fSNO}
\label{table:fSNO elliptic ub}
\resizebox{\textwidth}{!}{
\begin{tabular}{llrrrrrrrrrr}
\toprule
 equation & $N_{\text{train}}$  & $E_{\text{train}}$ & $E_{\text{test}}$ & $E^{\text{ub}}_{\text{train}}$ & $E^{\text{ub}}_{\text{test}}$ & $R^{\text{ub}}_{\text{train}}$ & $R^{\text{ub}}_{\text{test}}$ & $\widetilde{E}^{\text{ub}}_{\text{train}}$ & $\widetilde{E}^{\text{ub}}_{\text{test}}$ & $\widetilde{R}^{\text{ub}}_{\text{train}}$ & $\widetilde{R}^{\text{ub}}_{\text{test}}$ \\
\midrule
\multirow[t]{6}{*}{\eqref{eq:elliptic_1}, $J=5$} & 200 & 0.024 & 0.032 & 0.045 & 0.065 & 0.900 & 0.888 & 0.144 & 0.172 & 0.813 & 0.831 \\
 & 400 & 0.014 & 0.018 & 0.025 & 0.035 & 0.823 & 0.869 & 0.109 & 0.132 & 0.753 & 0.787 \\
 & 600 & 0.011 & 0.014 & 0.018 & 0.026 & 0.818 & 0.853 & 0.095 & 0.112 & 0.705 & 0.774 \\
 & 800 & 0.010 & 0.012 & 0.016 & 0.022 & 0.819 & 0.853 & 0.094 & 0.104 & 0.719 & 0.774 \\
 & 1000 & 0.009 & 0.011 & 0.015 & 0.020 & 0.824 & 0.849 & 0.089 & 0.100 & 0.722 & 0.756 \\
 & 1200 & 0.008 & 0.010 & 0.015 & 0.018 & 0.828 & 0.832 & 0.087 & 0.094 & 0.722 & 0.721 \\
\cline{1-12}
\multirow[t]{6}{*}{\eqref{eq:elliptic_1}, $J=6$} & 200 & 0.050 & 0.055 & 0.123 & 0.109 & 0.872 & 0.795 & 0.371 & 0.355 & 0.808 & 0.707 \\
\bottomrule
\end{tabular}
}
\end{table}

\begin{table}[!ht]
\centering
\caption{MLP}
\label{table:MLP elliptic ub}
\resizebox{\textwidth}{!}{
\begin{tabular}{llrrrrrrrrrr}
\toprule
 equation & $N_{\text{train}}$  & $E_{\text{train}}$ & $E_{\text{test}}$ & $E^{\text{ub}}_{\text{train}}$ & $E^{\text{ub}}_{\text{test}}$ & $R^{\text{ub}}_{\text{train}}$ & $R^{\text{ub}}_{\text{test}}$ & $\widetilde{E}^{\text{ub}}_{\text{train}}$ & $\widetilde{E}^{\text{ub}}_{\text{test}}$ & $\widetilde{R}^{\text{ub}}_{\text{train}}$ & $\widetilde{R}^{\text{ub}}_{\text{test}}$ \\
\midrule
\multirow[t]{6}{*}{\eqref{eq:elliptic_1}, $J=5$} & 200 & 0.037 & 0.037 & 0.069 & 0.070 & 0.878 & 0.826 & 0.377 & 0.353 & 0.854 & 0.847 \\
 & 400 & 0.028 & 0.026 & 0.048 & 0.051 & 0.889 & 0.846 & 0.207 & 0.204 & 0.744 & 0.720 \\
 & 600 & 0.022 & 0.021 & 0.036 & 0.039 & 0.862 & 0.823 & 0.201 & 0.199 & 0.675 & 0.684 \\
 & 800 & 0.019 & 0.018 & 0.031 & 0.032 & 0.835 & 0.785 & 0.187 & 0.187 & 0.553 & 0.513 \\
 & 1000 & 0.017 & 0.017 & 0.029 & 0.030 & 0.838 & 0.785 & 0.193 & 0.194 & 0.618 & 0.587 \\
 & 1200 & 0.016 & 0.016 & 0.027 & 0.027 & 0.861 & 0.814 & 0.238 & 0.231 & 0.875 & 0.867 \\
\cline{1-12}
\multirow[t]{6}{*}{\eqref{eq:elliptic_1}, $J=6$} & 200 & 0.043 & 0.044 & 0.076 & 0.072 & 0.852 & 0.843 & 0.368 & 0.363 & 0.816 & 0.796 \\
\bottomrule
\end{tabular}
}
\end{table}

\begin{table}[!ht]
\centering
\caption{UNet}
\label{table:UNet elliptic ub}
\resizebox{\textwidth}{!}{
\begin{tabular}{llrrrrrrrrrr}
\toprule
 equation & $N_{\text{train}}$  & $E_{\text{train}}$ & $E_{\text{test}}$ & $E^{\text{ub}}_{\text{train}}$ & $E^{\text{ub}}_{\text{test}}$ & $R^{\text{ub}}_{\text{train}}$ & $R^{\text{ub}}_{\text{test}}$ & $\widetilde{E}^{\text{ub}}_{\text{train}}$ & $\widetilde{E}^{\text{ub}}_{\text{test}}$ & $\widetilde{R}^{\text{ub}}_{\text{train}}$ & $\widetilde{R}^{\text{ub}}_{\text{test}}$ \\
\midrule
\multirow[t]{6}{*}{\eqref{eq:elliptic_1}, $J=5$} & 200 & 0.020 & 0.026 & 0.040 & 0.061 & 0.889 & 0.868 & 0.271 & 0.279 & 0.833 & 0.861 \\
 & 400 & 0.015 & 0.020 & 0.027 & 0.043 & 0.875 & 0.887 & 0.177 & 0.185 & 0.860 & 0.840 \\
 & 600 & 0.010 & 0.016 & 0.019 & 0.032 & 0.850 & 0.881 & 0.151 & 0.163 & 0.812 & 0.809 \\
 & 800 & 0.009 & 0.014 & 0.016 & 0.027 & 0.865 & 0.876 & 0.103 & 0.111 & 0.830 & 0.798 \\
 & 1000 & 0.008 & 0.013 & 0.015 & 0.025 & 0.859 & 0.874 & 0.095 & 0.104 & 0.865 & 0.813 \\
 & 1200 & 0.007 & 0.011 & 0.012 & 0.021 & 0.869 & 0.866 & 0.086 & 0.094 & 0.813 & 0.754 \\
\cline{1-12}
\multirow[t]{6}{*}{\eqref{eq:elliptic_1}, $J=6$} & 200 & 0.033 & 0.037 & 0.192 & 0.136 & 0.900 & 0.802 & 0.535 & 0.496 & 0.800 & 0.665 \\
\bottomrule
\end{tabular}
}
\end{table}

\begin{table}[!h]
  \caption{Results for 1D unsupervised learning equations averaged with respect to architectures, $N_{\text{train}}=1200$.}
  \label{table:unsup_1D_eq}
  \centering
\begin{tabular}{crrrrrr}
\toprule
{equation} &  $E_{\text{train}}$ &  $E_{\text{test}}$ &  $E^{\text{ub}}_{\text{train}}$ &  $E^{\text{ub}}_{\text{test}}$ &  $R^{\text{ub}}_{\text{train}}$ &  $R^{\text{ub}}_{\text{test}}$ \\
\midrule
\eqref{eq:elliptic_1_1D} &               0.011 &              0.015 &                           0.019 &                          0.034 &                           0.827 &                          0.855 \\
\eqref{eq:elliptic_2_1D} &               0.010 &              0.015 &                           0.015 &                          0.032 &                           0.822 &                          0.854 \\
\eqref{eq:elliptic_3_1D} &               0.033 &              0.034 &                           0.037 &                          0.039 &                           0.791 &                          0.792 \\
\bottomrule
\end{tabular}
\end{table}
\begin{table}[!th]
  \caption{Results for 1D unsupervised learning equations averaged with respect to equations, $N_{\text{train}}=1200$.}
  \label{table:unsup_1D_m}
  \centering
\begin{tabular}{crrrrrr}
\toprule
Architecture &  $E_{\text{train}}$ &  $E_{\text{test}}$ &  $E^{\text{ub}}_{\text{train}}$ &  $E^{\text{ub}}_{\text{test}}$ &  $R^{\text{ub}}_{\text{train}}$ &  $R^{\text{ub}}_{\text{test}}$ \\
\midrule
fSNO      &               0.014 &              0.016 &                           0.022 &                          0.026 &                           0.828 &                          0.861 \\
ChebNO    &               0.016 &              0.018 &                           0.026 &                          0.029 &                           0.835 &                          0.853 \\
DilResNet &               0.017 &              0.018 &                           0.018 &                          0.022 &                           0.771 &                          0.777 \\
UNet      &               0.023 &              0.033 &                           0.028 &                          0.061 &                           0.821 &                          0.844 \\
\bottomrule
\end{tabular}
\end{table}
\begin{table}[!th]
  \caption{Results for 1D unsupervised learning for distinct $N_{\text{train}}$  averaged with respect to the network type and equations.}
  \label{table:unsup_1D_trsize}
  \centering
\begin{tabular}{crrrrrr}
\toprule
$N_{\text{train}}$ &  $E_{\text{train}}$ &  $E_{\text{test}}$ &  $E^{\text{ub}}_{\text{train}}$ &  $E^{\text{ub}}_{\text{test}}$ &  $R^{\text{ub}}_{\text{train}}$ &  $R^{\text{ub}}_{\text{test}}$ \\
\midrule
200      &               0.036 &              0.048 &                           0.060 &                          0.102 &                           0.840 &                          0.839 \\
400      &               0.026 &              0.033 &                           0.041 &                          0.066 &                           0.825 &                          0.841 \\
600      &               0.023 &              0.029 &                           0.032 &                          0.052 &                           0.810 &                          0.839 \\
800      &               0.021 &              0.026 &                           0.030 &                          0.047 &                           0.823 &                          0.842 \\
1000     &               0.019 &              0.023 &                           0.026 &                          0.040 &                           0.823 &                          0.845 \\
1200     &               0.018 &              0.021 &                           0.023 &                          0.035 &                           0.814 &                          0.834 \\
\bottomrule
\end{tabular}
\end{table}
\begin{table}[!th]
  \caption{Results for 1D unsupervised learning for distinct $N_{\text{train}}$  averaged with respect to equations.}
  \label{table:unsup_1D_unet}
  \centering
\begin{tabular}{cllllllll}
\toprule
$N_{\text{train}}$  &                \multicolumn{2}{c}{fSNO} & \multicolumn{2}{c}{ChebNO} & \multicolumn{2}{c}{DilResNet} & \multicolumn{2}{c}{UNet} \\
{} &  $E_{\text{train}}$ &     $E_{\text{test}}$ &  $E_{\text{train}}$ &     $E_{\text{test}}$ &          $E_{\text{train}}$ &     $E_{\text{test}}$ &                        $E_{\text{train}}$ &    $E_{\text{test}}$ \\
\midrule
200                &               0.028 &  0.036 &                       0.035 &                 0.044 &                                      0.031 &                 0.041 &                        0.05 &  0.071 \\
400                &                0.02 &                 0.025 &                       0.025 &                 0.029 &                                       0.02 &                 0.024 &                       0.038 &                0.055 \\
600                &               0.019 &                 0.022 &                        0.02 &                 0.023 &                                      0.019 &                 0.021 &                       0.033 &                0.049 \\
800                &               0.017 &                 0.019 &                       0.019 &                 0.021 &                       0.018 &                  0.02 &                       0.031 &                0.044 \\
1000               &               0.015 &                 0.016 &                       0.017 &  0.018 &                                      0.017 &                 0.019 &                       0.028 &                 0.04 \\
1200               &               0.014 &                 0.016 &                       0.016 &  0.018 &                                      0.017 &  0.018 &                       0.023 &                0.033 \\
\bottomrule
\end{tabular}
\end{table}

\begin{table}[!ht]
\centering
\caption{fSNO train and test errors with different  $\lambda$  for 1D datasets}
\label{table:fsno_lambda}
\begin{tabular}{lllllllllll}
\toprule
equation &              \multicolumn{2}{c}{$\lambda = 0.01$} &                                \multicolumn{2}{c}{$\lambda = 0.1$} &                                \multicolumn{2}{c}{$\lambda = 1$} &                                   \multicolumn{2}{c}{$\lambda = 10$} &                             \multicolumn{2}{c}{$\lambda = 100$}                \\
                   &  $E_{\text{train}}$ &     $E_{\text{test}}$ &    $E_{\text{train}}$ &     $E_{\text{test}}$ &  $E_{\text{train}}$ &     $E_{\text{test}}$ &    $E_{\text{train}}$ &  $E_{\text{test}}$ &  $E_{\text{train}}$ &     $E_{\text{test}}$ \\
\midrule
\eqref{eq:elliptic_1_1D}  &                 --- &                --- &  0.009 &              0.011 &  0.009 &              0.011 &               0.011 &              0.012 &               0.015 &                 0.017 \\
\eqref{eq:elliptic_2_1D}  &               0.028 &               0.03 &  0.009 &              0.011 &  0.009 &              0.011 &                0.01 &              0.011 &               0.016 &  0.018 \\
\eqref{eq:elliptic_3_1D}  &               0.058 &              0.061 &                 0.042 &              0.044 &                  0.02 &              0.021 &               0.023 &              0.024 &               0.082 &                 0.083 \\
\bottomrule
\end{tabular}
\end{table}

\begin{table}[!ht]
\centering
\caption{ChebNO train and test errors with different  $\lambda$  for 1D datasets}
\label{table:chebno_lambda}
\begin{tabular}{lllllllllll}
\toprule
equation &              \multicolumn{2}{c}{$\lambda = 0.01$} &                                \multicolumn{2}{c}{$\lambda = 0.1$} &                                \multicolumn{2}{c}{$\lambda = 1$} &                                   \multicolumn{2}{c}{$\lambda = 10$} &                             \multicolumn{2}{c}{$\lambda = 100$}                \\
                   &  $E_{\text{train}}$ &     $E_{\text{test}}$ &    $E_{\text{train}}$ &     $E_{\text{test}}$ &  $E_{\text{train}}$ &     $E_{\text{test}}$ &    $E_{\text{train}}$ &  $E_{\text{test}}$ &  $E_{\text{train}}$ &     $E_{\text{test}}$ \\
\midrule
\eqref{eq:elliptic_1_1D} &                 --- &              0.058 &                0.02 &               0.02 &                 0.024 &              0.023 &               0.029 &              0.029 &               0.148 &                0.147 \\
\eqref{eq:elliptic_2_1D} &               0.055 &              0.056 &               0.014 &              0.014 &                 0.015 &              0.015 &               0.031 &              0.033 &               0.133 &  0.142 \\
\eqref{eq:elliptic_3_1D} &               0.078 &              0.081 &               0.053 &              0.056 &  0.0360 &              0.037 &               0.065 &              0.065 &               0.109 &                0.109 \\
\bottomrule
\end{tabular}
\end{table}

\begin{table}[!ht]
\centering
\caption{DilResNet train and test errors with different $\lambda$ for 1D datasets}
\label{table:Dilresnet_lambda}
\begin{tabular}{lllllllllll}
\toprule
equation &              \multicolumn{2}{c}{$\lambda = 0.01$} &                              \multicolumn{2}{c}{$\lambda = 0.1$} &                              \multicolumn{2}{c}{$\lambda = 1$} &                                 \multicolumn{2}{c}{$\lambda = 10$} &                           \multicolumn{2}{c}{$\lambda = 100$}                \\
                   &  $E_{\text{train}}$ &     $E_{\text{test}}$ &    $E_{\text{train}}$ &     $E_{\text{test}}$ &  $E_{\text{train}}$ &     $E_{\text{test}}$ &    $E_{\text{train}}$ &  $E_{\text{test}}$ &  $E_{\text{train}}$ &     $E_{\text{test}}$ \\
\midrule
\eqref{eq:elliptic_1_1D}  &                 --- &                   --- &                 0.008 &  0.009 &               0.008 &                  0.01 &  0.009 &               0.01 &               0.012 &  0.013 \\
\eqref{eq:elliptic_2_1D}  &               0.012 &  0.013 &  0.009 &                  0.01 &               0.008 &  0.009 &                 0.008 &               0.01 &               0.016 &  0.018 \\
\eqref{eq:elliptic_3_1D}  &               0.085 &   0.086 &   0.071 &                 0.073 &               0.035 &  0.036 &                 0.048 &              0.047 &               0.048 &                 0.047 \\
\bottomrule
\end{tabular}
\end{table}

\begin{table}[!ht]
\centering
\caption{UNet train and test errors with different  $\lambda$ for 1D datasets}
\label{table:u_lambda}
\begin{tabular}{lllllllllll}
\toprule
equation &              \multicolumn{2}{c}{$\lambda = 0.01$}                &                 \multicolumn{2}{c}{$\lambda = 0.1$}                &                 \multicolumn{2}{c}{$\lambda = 1$}                  &                  \multicolumn{2}{c}{$\lambda = 10$}               &               \multicolumn{2}{c}{$\lambda = 100$}                \\
                   &  $E_{\text{train}}$ &     $E_{\text{test}}$ &    $E_{\text{train}}$ &     $E_{\text{test}}$ &  $E_{\text{train}}$ &     $E_{\text{test}}$ &    $E_{\text{train}}$ &  $E_{\text{test}}$ &  $E_{\text{train}}$ &     $E_{\text{test}}$ \\
\midrule
\eqref{eq:elliptic_1_1D}  &  0.026 &                 0.028 &               0.167 &              0.164 &                 0.034 &              0.039 &               0.044 &              0.048 &               0.168 &              0.166 \\
\eqref{eq:elliptic_2_1D}  &                 0.021 &  0.036 &               0.032 &              0.048 &  0.026 &              0.033 &               0.038 &              0.047 &                0.16 &              0.171 \\
\eqref{eq:elliptic_3_1D}  &                 0.114 &                 0.115 &                0.16 &              0.161 &                 0.062 &              0.062 &                0.11 &               0.11 &               0.138 &              0.139 \\
\bottomrule
\end{tabular}
\end{table}

\begin{table}[!ht]
\centering
\caption{fSNO  train and test results for all 1D datasets. }
\label{table:fsno_datasets}
\begin{tabular}{lrrrrrr}
\toprule
equation  &  $E_{\text{train}}$ &  $E_{\text{test}}$ &  $E^{\text{up}}_{\text{train}}$ &  $E^{\text{up}}_{\text{test}}$ &  $R^{\text{up}}_{\text{train}}$ &  $R^{\text{up}}_{\text{test}}$ \\
\midrule
\eqref{eq:elliptic_1_1D} &               0.009 &              0.011 &                           0.017 &                          0.023 &                           0.829 &                          0.869 \\
\eqref{eq:elliptic_2_1D} &               0.009 &              0.011 &                           0.016 &                          0.020 &                           0.819 &                          0.846 \\
\eqref{eq:elliptic_3_1D} &               0.023 &              0.025 &                           0.033 &                          0.035 &                           0.848 &                          0.879 \\
\bottomrule
\end{tabular}
\end{table}

\begin{table}[!ht]
\centering
\caption{ChebNO  train and test results for all 1D datasets. }
\label{table:chebno_datasets}
\begin{tabular}{lrrrrrr}
\toprule
equation &  $E_{\text{train}}$ &  $E_{\text{test}}$ &  $E^{\text{up}}_{\text{train}}$ &  $E^ {\text{up}}_{\text{test}}$ &  $R^ {\text{up}}_{\text{train}}$ &  $R^ {\text{up}}_{\text{test}}$ \\
\midrule
\eqref{eq:elliptic_1_1D} &               0.024 &              0.023 &                           0.051 &                          0.050 &                           0.826 &                          0.833 \\
\eqref{eq:elliptic_2_1D} &               0.015 &              0.015 &                           0.033 &                          0.034 &                           0.812 &                          0.817 \\
\eqref{eq:elliptic_3_1D} &               0.036 &              0.037 &                           0.049 &                          0.051 &                           0.896 &                          0.914 \\
\bottomrule
\end{tabular}
\end{table}

\begin{table}[!ht]
\centering
\caption{DilResNet train and test results for all 1D datasets.}
\label{table:Dilresnet_datasets}
\begin{tabular}{lrrrrrr}
\toprule
equation &  $E_{\text{train}}$ &  $E_{\text{test}}$ &  $E^{\text{up}}_{\text{train}}$ &  $E^{\text{up}}_{\text{test}}$ &  $R^{\text{up}}_{\text{train}}$ &  $R^{\text{up}}_{\text{test}}$ \\

\midrule
\eqref{eq:elliptic_1_1D} &               0.008 &              0.009 &                           0.012 &                          0.016 &                           0.795 &                          0.825 \\
\eqref{eq:elliptic_2_1D} &               0.008 &              0.009 &                           0.012 &                          0.016 &                           0.825 &                          0.847 \\
\eqref{eq:elliptic_3_1D} &               0.035 &              0.036 &                           0.031 &                          0.033 &                           0.693 &                          0.658 \\
\bottomrule
\end{tabular}
\end{table}

\begin{table}[!ht]
\centering
\caption{UNet  train and test results for all 1D datasets.}
\label{table:u_datasets}
\begin{tabular}{lrrrrrr}
\toprule
equation &  $E_{\text{train}}$ &  $E_{\text{test}}$ &  $E^{\text{up}}_{\text{train}}$ &  $E^{\text{up}}_{\text{test}}$ &  $R^{\text{up}}_{\text{train}}$ &  $R^{\text{up}}_{\text{test}}$ \\

\midrule
\eqref{eq:elliptic_1_1D} &               0.013 &              0.026 &                           0.021 &                          0.067 &                           0.846 &                          0.894 \\
\eqref{eq:elliptic_2_1D} &               0.011 &              0.027 &                           0.016 &                          0.069 &                           0.847 &                          0.896 \\
\eqref{eq:elliptic_3_1D} &               0.045 &              0.045 &                           0.045 &                          0.047 &                           0.769 &                          0.743 \\
\bottomrule
\end{tabular}
\end{table}

\begin{table}[!ht]
\centering
\caption{fSNO train and test errors with different  $\lambda$  for 2D datasets}
\label{table:fsno_lambda_2D}
\begin{tabular}{lllllllllll}
\toprule
equation &              \multicolumn{2}{c}{$\lambda = 0.01$} &                                \multicolumn{2}{c}{$\lambda = 0.1$} &                                \multicolumn{2}{c}{$\lambda = 1$} &                                   \multicolumn{2}{c}{$\lambda = 10$} &                             \multicolumn{2}{c}{$\lambda = 100$}                \\
                   &  $E_{\text{train}}$ &     $E_{\text{test}}$ &    $E_{\text{train}}$ &     $E_{\text{test}}$ &  $E_{\text{train}}$ &     $E_{\text{test}}$ &    $E_{\text{train}}$ &  $E_{\text{test}}$ &  $E_{\text{train}}$ &     $E_{\text{test}}$ \\
\midrule
\eqref{eq:elliptic_1} &  0.026 &               0.027 &                0.016 &  0.018 &                0.014 &                 0.016 &               0.017 &               0.019 &  0.104 &  0.102 \\
\eqref{eq:elliptic_3} &                 0.028 &               0.028 &                 0.02 &                  0.02 &                0.019 &                 0.019 &               0.021 &               0.021 &                 0.03 &                0.029 \\
\eqref{eq:elliptic_4} &                 0.021 &               0.023 &                0.015 &                 0.017 &                0.016 &  0.018 &               0.017 &               0.019 &  0.071 &                 0.070 \\
\eqref{eq:elliptic_6} &                 0.032 &               0.033 &                0.016 &                 0.019 &                0.014 &                 0.017 &               0.017 &               0.019 &                0.105 &                0.105 \\
\bottomrule
\end{tabular}
\end{table}

\begin{table}[!ht]
\centering
\caption{ChebNO train and test errors with different  $\lambda$ for 2D datasets}
\label{table:chebno_lambda_2D}
\begin{tabular}{lllllllllll}
\toprule
equation &              \multicolumn{2}{c}{$\lambda = 0.01$} &                                \multicolumn{2}{c}{$\lambda = 0.1$} &                                \multicolumn{2}{c}{$\lambda = 1$} &                                   \multicolumn{2}{c}{$\lambda = 10$} &                             \multicolumn{2}{c}{$\lambda = 100$}                \\
                   &  $E_{\text{train}}$ &     $E_{\text{test}}$ &    $E_{\text{train}}$ &     $E_{\text{test}}$ &  $E_{\text{train}}$ &     $E_{\text{test}}$ &    $E_{\text{train}}$ &  $E_{\text{test}}$ &  $E_{\text{train}}$ &     $E_{\text{test}}$ \\
\midrule
\eqref{eq:elliptic_1} &  0.026 &                 0.028 &                0.016 &  0.018 &                0.015 &                 0.017 &               0.016 &  0.018 &                0.064 &                 0.062 \\
\eqref{eq:elliptic_3} &                 0.032 &                 0.031 &                0.021 &                 0.021 &                 0.02 &                  0.02 &               0.021 &                 0.021 &                0.028 &                 0.027 \\
\eqref{eq:elliptic_4} &                 0.024 &  0.026 &                0.016 &  0.018 &                0.016 &  0.018 &                0.02 &                 0.021 &                0.055 &                 0.054 \\
\eqref{eq:elliptic_6} &                  0.03 &                 0.031 &                0.017 &                 0.019 &                0.016 &                 0.017 &               0.019 &                 0.021 &                0.025 &  0.026 \\
\bottomrule
\end{tabular}
\end{table}

\begin{table}[!ht]
\centering
\caption{DilResNet train and test errors with different  $\lambda$ for 2D datasets}
\label{table:dilresnet_lambda_2D}
\begin{tabular}{lllllllllll}
\toprule
equation &              \multicolumn{2}{c}{$\lambda = 0.01$} &                                \multicolumn{2}{c}{$\lambda = 0.1$} &                                \multicolumn{2}{c}{$\lambda = 1$} &                                   \multicolumn{2}{c}{$\lambda = 10$} &                             \multicolumn{2}{c}{$\lambda = 100$}                \\
                   &  $E_{\text{train}}$ &     $E_{\text{test}}$ &    $E_{\text{train}}$ &     $E_{\text{test}}$ &  $E_{\text{train}}$ &     $E_{\text{test}}$ &    $E_{\text{train}}$ &  $E_{\text{test}}$ &  $E_{\text{train}}$ &     $E_{\text{test}}$ \\
\midrule
\eqref{eq:elliptic_1} &  0.013 &  0.013 &                  0.01 &                  0.01 &  0.009 &                0.01 &  0.009 &                  0.01 &                 0.017 &                 0.017 \\
\eqref{eq:elliptic_3} &                 0.021 &                  0.02 &  0.018 &  0.018 &  0.018 &               0.018 &  0.018 &  0.018 &                 0.027 &  0.026 \\
\eqref{eq:elliptic_4} &  0.013 &                 0.014 &                 0.012 &                 0.012 &                 0.011 &               0.012 &                 0.012 &                 0.012 &                 0.048 &                 0.047 \\
\eqref{eq:elliptic_6} &  0.013 &                 0.014 &  0.009 &                  0.01 &                  0.01 &                0.01 &                  0.01 &                 0.011 &  0.026 &  0.026 \\
\bottomrule
\end{tabular}
\end{table}

\begin{table}[!ht]
\centering
\caption{UNet train and test errors with different  $\lambda$ for 2D datasets}
\label{table:unet_lambda_2D}
\begin{tabular}{lllllllllll}
\toprule
equation &              \multicolumn{2}{c}{$\lambda = 0.01$} &                                \multicolumn{2}{c}{$\lambda = 0.1$} &                                \multicolumn{2}{c}{$\lambda = 1$} &                                   \multicolumn{2}{c}{$\lambda = 10$} &                             \multicolumn{2}{c}{$\lambda = 100$}                \\
                   &  $E_{\text{train}}$ &     $E_{\text{test}}$ &    $E_{\text{train}}$ &     $E_{\text{test}}$ &  $E_{\text{train}}$ &     $E_{\text{test}}$ &    $E_{\text{train}}$ &  $E_{\text{test}}$ &  $E_{\text{train}}$ &     $E_{\text{test}}$ \\
\midrule
\eqref{eq:elliptic_1} &                 0.02 &               0.021 &  0.013 &                 0.015 &                0.015 &               0.015 &                0.021 &               0.022 &  0.104 &  0.102 \\
\eqref{eq:elliptic_3} &                 0.04 &                0.04 &                 0.028 &                 0.027 &                0.019 &               0.019 &                0.091 &                0.09 &                 0.11 &                0.109 \\
\eqref{eq:elliptic_4} &                0.024 &               0.024 &   0.071 &                  0.07 &  0.071 &                0.07 &  0.071 &                0.07 &  0.071 &                 0.07 \\
\eqref{eq:elliptic_6} &                0.016 &               0.017 &                 0.016 &  0.018 &                0.012 &               0.015 &                0.032 &               0.035 &                0.105 &                0.105 \\
\bottomrule
\end{tabular}
\end{table}

\begin{table}[!ht]
\centering
\caption{fSNO  train and test results for all 2D datasets. }
\label{table:fsno_datasets_2d}
\begin{tabular}{lrrrrrr}
\toprule
equation  &  $E_{\text{train}}$ &  $E_{\text{test}}$ &  $E^{\text{up}}_{\text{train}}$ &  $E^{\text{up}}_{\text{test}}$ &  $R^{\text{up}}_{\text{train}}$ &  $R^{\text{up}}_{\text{test}}$ \\
\midrule
\eqref{eq:elliptic_1} &                0.014 &               0.016 &                            0.029 &                           0.041 &                            0.840 &                           0.790 \\
\eqref{eq:elliptic_3} &                0.019 &               0.019 &                            0.026 &                           0.027 &                            0.829 &                           0.821 \\
\eqref{eq:elliptic_4} &                0.016 &               0.017 &                            0.025 &                           0.031 &                            0.946 &                           0.942 \\
\eqref{eq:elliptic_6} &                0.014 &               0.017 &                            0.031 &                           0.042 &                            0.862 &                           0.825 \\
\bottomrule
\end{tabular}
\end{table}

\begin{table}[!ht]
\centering
\caption{ChebNO  train and test results for all 2D datasets. }
\label{table:chebno_datasets_2d}
%
\end{table}

\begin{table}[!ht]
\centering
\caption{DilResNet  train and test results for all 2D datasets. }
\label{table:dilresnet_datasets_2d}
%
\end{table}

\begin{table}[!ht]
\centering
\caption{UNet  train and test results for all 2D datasets. }
\label{table:unet_datasets_2d}
%
\end{table}

\bibliography{refs.bib}